\documentclass[twocolumn]{autart} 
\usepackage{graphicx} 
\usepackage{color,natbib} 
\usepackage{mdframed,subfigure}

\usepackage{amsmath,amssymb,mathrsfs}
\usepackage{mathtools,xparse}

\usepackage[dvipsnames]{xcolor}

\newcommand{\PPP}{{\mathbf P}}

\DeclarePairedDelimiter{\norm}{\lVert}{\rVert}
\newcommand{\rfb}[1]{\mbox{\rm
(\ref{#1})}\ifx\undefined\stillediting\else:\fbox{$#1$}\fi} 
\newcommand{\rarrow} {\mathop{\rightarrow}}                
\newcommand{\FORALL} {{\hbox{$\hskip 11mm \forall \;$}}}
\renewcommand{\e}    {{\varepsilon}}
\newcommand{\m}      {{\hbox{\hskip 1pt}}}
\newcommand{\mm}     {{\hbox{\hskip 0.5pt}}}
\newcommand{\nm}     {{\hbox{\hskip -3pt}}}
\newcommand{\sbluff} {{\hbox{\raise 10pt \hbox{\mm}}}}
\renewcommand{\l}    {{\lambda}}
\renewcommand{\L}    {{\Lambda}}
\renewcommand{\Re}   {{\rm Re\,}}

\newcommand{\rline}  {{\mathbb R}}
\newcommand{\dd}     {{\rm d}\hbox{\hskip 0.5pt}}

\renewcommand{\o}  {{\omega}} 

\newcommand{\Dscr} {{\mathcal D}}

\newcommand{\Nscr} {{\mathcal N}}
\newcommand{\Uscr} {{\mathcal U}} 
\newcommand{\Vscr} {{\mathcal V}} 
\newcommand{\Rscr} {{\mathcal R}}

\newcommand{\Om}   {{\Omega}}

\DeclareMathOperator*{\argmin}{arg\,min}
\DeclareMathOperator*{\argmax}{arg\,max}
\newcommand{\bbm}[1]{\left[\begin{matrix} #1 
   \end{matrix}\right]}
\newcommand{\sbm}[1]{\left[\begin{smallmatrix} #1 
   \end{smallmatrix}\right]}

\numberwithin{equation}{section}

\newtheorem{theorem}{Theorem}[section]

\newtheorem{lemma}[theorem]{Lemma}
\newtheorem{proposition}[theorem]{Proposition}
\newtheorem{definition}[theorem]{Definition}

\newtheorem{remark}[theorem]{Remark}

\begin{document}
\begin{frontmatter}
\title{PI control of stable nonlinear plants using\\ projected
   dynamical systems\thanksref{footnoteinfo}} 

\thanks[footnoteinfo]{The authors are team-members in the 
ITN network ConFlex. This project has received funding from the 
European Union's Horizon 2020 research and innovation program 
under the Marie Sklodowska-Curie grant agreement no. 765579. The 
research was also supported by grant no. 2802/21 from the
Israel Science Foundation (ISF).}

\author[TAU]{Pietro Lorenzetti}\ead{plorenzetti@tauex.tau.ac.il},
\author[TAU]{George Weiss}\ead{gweiss@tauex.tau.ac.il},
\address[TAU]{School of Electrical Engineering, Tel Aviv 
   University, Ramat Aviv 69978, Israel }  


\begin{abstract} 
This paper presents a novel anti-windup proportional-integral
controller for stable multi-input multi-output nonlinear plants. We
use tools from projected dynamical systems theory to force the
integrator state to remain in a desired (compact and convex) region,
such that the plant input steady-state values satisfy the operational
constraints of the problem. Under suitable monotonicity
assumptions on the plant steady-state input-output map, we use
singular perturbation theory results to prove the existence of a
sufficiently small controller gain ensuring closed-loop (local)
exponential stability and reference tracking for a feasible set of
constant references.  We suggest a particular controller design, which
embeds (when possible) the right inverse of the plant steady-state
input-output map. The relevance of the proposed controller scheme is
validated through an application in the power systems domain, namely,
the output (active and reactive) power regulation for a grid-connected
synchronverter.
\end{abstract}

\end{frontmatter}

\section{Introduction} \label{sec:1} 

One of the fundamental problems in control theory is the {\em
regulator problem}, where the objective is to design a controller that
forces the output of a plant to track a reference signal, while
rejecting possible disturbances.  It is often convenient to assume
that the reference and the disturbance signals are generated by a
fictitious system called the {\em exosystem}, which is an expression
of our previous knowledge about these signals. When it is reasonable
to assume that the signals originate from an exosystem, then the {\em
internal model principle} (see \cite{Davison1976,Francis1975}), in its
version for linear time-invariant (LTI) systems, states that the
regulator problem is solved if all the unstable eigenvalues of the
exosystem are poles of the controller. In the case of constant
signals, this principle suggests that an integral controller is needed
to solve the regulator problem.

For LTI systems, if the plant DC-gain sign is known and the plant is
stable, then, for sufficiently small controller gains of suitable
sign, the closed-loop system formed by the plant with the integral
controller is stable, and the regulator problem is solved (see
\cite{Morari1985}). Similar results have been established in
\cite{Desoer1985} for multi-input multi-output (MIMO) globally stable
nonlinear systems, using singular perturbation (SP) theory. In their work,
the plant DC-gain sign assumption is replaced with an assumption on
the monotonicity of the plant steady-state input-output mapping. The
same approach has been used in our paper \cite{Lor2020} for the
stability analysis of a stable nonlinear plant connected in feedback
with a single-input single-output (SISO) low-gain anti-windup (AW)
proportional-integral (PI) controller. The recent paper
\cite{Simpson-Porco2020} on low-gain integral control for stable nonlinear
systems, also employing SP tools, has
generalized the main result from \cite{Desoer1985}. The assumption on
the monotonicity of the plant steady-state input-output mapping has
been replaced with the uniform infinitesimal contracting property of
the reduced dynamics. It has then been shown how this relaxed
assumption recovers the one in \cite{Morari1985} for a linear plant.
(For an extension of this result to nonlinear discrete-time systems
with input constraints see \cite{Simpson-Porco2021}.)  Other
interesting results for linear systems with input-output
nonlinearities and low-gain integral controllers are in
\cite{Logemann1999,Guiver2017}.

An ubiquitous problem in control applications is \textit{windup}. This
happens when there is a mismatch between the controller output and the
actual plant input, e.g., due to actuator limitations, causing long
transients, oscillations and even instability (\cite{Kothare1994}).
Several AW techniques have been proposed, resulting in a vast
literature on the topic. It is not in the scope of this paper to
provide a detailed review on AW control, instead we refer to
\cite{Astrom1989,Edwards1998,Kothare1994,Tarbouriech2009,
Zaccarian2002} (and the references therein). Even though such AW control strategies prove to
be effective, they mainly deal with linear plants and they often 
require the solution of an LMI-based optimization problem.

Our aim is to formulate a simple novel AW MIMO low-gain PI
controller for stable nonlinear systems, which forces the integrator
state (and, thus, the steady-state plant input) to stay in a
desired (compact and convex) region, using tools from projected
dynamical systems (PDS) theory. We provide a rigorous closed-loop
stability analysis, following the approach of
\cite{Desoer1985,Lor2020,Simpson-Porco2020}, by using SP
theory results. We derive a sufficient condition on the
controller gain that ensures closed-loop stability and reference
tracking for a ``naturally'' feasible set of constant references. This
result generalizes \cite{Lor2020}, where a similar PI AW control
strategy, using the saturating integrator introduced in
\cite{Lor2020Conf}, was formulated for SISO stable nonlinear plants.
The PI SISO saturating integrator performs remarkably well in several
applications, see, e.g., \cite{Lor2020,Lor2021,Natarajan2017}. We expect similar successful
performance (and more) for its MIMO formulation shown in Fig.~\ref{fig:clos_MIMO}. A preliminary
version of this paper, considering only integral control 
$(\tau_p=0)$, without the
block $\Nscr$, and with no numerical example,
has been presented in
our recent conference paper \cite{Lor2021Conf}. Here, we introduce an
additional degree of freedom in the controller design, namely, the
block $\mathcal{N}$ from Fig.~\ref{fig:clos_MIMO}, which we exploit to
embed (when possible) the right inverse of the plant steady-state
input-output map in the controller. We use a numerical example from the
power electronics domain to illustrate the relevance of the
proposed controller design.

Related to our work is the recent contribution \cite{Wang2020}, where
a bounded MIMO integral controller is presented (not based on PDS
theory), called bounded integral controller (BIC) (this is an
extension of the SISO BIC presented in \cite{Konstantopoulos2016}).
The MIMO BIC enforces a sum-of-squares-type of constraint with
time-varying input weights for the controller states, and
input-to-state practical stability (ISpS) is guaranteed (using a
small-gain argument) when the BIC is connected in feedback with an
ISpS nonlinear plant. On the other hand, the (natural) link between
PDS theory and AW design has been pointed out before. Investigations
in a similar direction have been carried out in \cite{Teo2011}, 
where a gradient projection AW (GPAW) scheme has been
proposed. In particular, they provide sufficient conditions under
which the following holds: Assume that the nominal unconstrained
closed-loop system is stable and achieves tracking in correspondence
of a certain equilibrium point, then the region of attraction of the
same equilibrium point in the associated GPAW closed-loop system is
``larger'' than the one of the unconstrained system. Although this
result is of interest, the sufficient conditions provided are
difficult to verify in practice (as pointed out by
the authors), see \cite[Theorem 2]{Teo2011}. Recently, the connection between PDS and AW schemes has
also been investigated in \cite{Hauswirth2020b}, with
applications to feedback optimization problems. In particular, they
show that the closed-loop solutions of a high-gain integral AW control
scheme uniformly converge to those of a PDS as the gain tends to
infinity, see Remark~\ref{rmk:Haws}. Finally, a low-gain projected integral control scheme for
exponentially stable discrete-time nonlinear systems is presented in
the recent contribution \cite{Simpson-Porco2021}.

The paper is organized as follows. In Sect.~\ref{sec2} we present some
background on PDS theory. In Sect.~\ref{sec3} the PI SISO
saturating integrator from \cite{Lor2020} is reformulated as a MIMO
controller, using tools from PDS theory, and the control problem 
is described in precise terms. Sect.~\ref{sec4} contains
our main result: the stability of the closed-loop
system (and the consequent reference tracking), proved using
SP theory. Finally, in Sect.~\ref{sec5} we
illustrate the performance of the proposed controller through an
application from power electronics.

\section{Background on PDS theory} \label{sec2} 

We present some background on PDS
theory, taken mainly from \cite[Ch.~2]{Nagurney1995}. \vspace{-1mm}

\textbf{Notation.} Let $X\subset\mathbb{R}^q$ be closed and
convex. Denote the boundary (interior) of $X$ by $\partial X$ (${\rm
int}\m X$). Define the set of \textit{inward normals} to $X$ at
$x\in\partial X$ by \vspace{-1mm}
\begin{equation*} \label{def:in_norm} 
   n(x) \m=\m \{\gamma \ \big|\ \norm{\gamma}=1,\m\m\text{and}\m\m\langle
\gamma,x-y\rangle\leq 0,\m\forall\m y \in X\}. \vspace{-1mm}
\end{equation*}

\begin{definition} \label{def:Pi} {\rm (\cite{Nagurney1995})}
Let $X\subset\rline^q$ be closed and convex, and let $z\in X$, 
$v\in\rline^q$. Define the projection operator $P_X$
onto $X$ as \vspace{-1mm}
\begin{equation}\label{eq:proj} 
   P_X(v) \m=\m \argmin_{w\in X} \norm{v-w}, \vspace{-1mm}
\end{equation}
and the directional derivative of $P_X$ at $z$, along $v$, as
\vspace{-1mm}
\begin{equation} \label{eq:Pi}
   \Pi_X(z,v) \m=\m \lim_{\delta\to 0^+} \frac{P_X(z+\delta v)-z}
{\delta}. \vspace{-1mm}
\end{equation}
\end{definition}

\vspace{-1mm} 
{\color{blue} \begin{lemma} \label{lmm:Pi}{\rm \cite[Lemma~2.1]{Nagurney1995}.} 
Let $X,z,v,\Pi$ be as in Definition \ref{def:Pi}. Then: \vspace{-3mm}
\begin{enumerate} 
\item If $z\in{\rm int}\m X$, then \ \ $\Pi_X(z,v)=v$. 
\item if $z\in\partial X$, then \ \ $\Pi_X(z,v)=v+\beta(z)n^*(z)$, where
\end{enumerate}
\vspace{-2mm}
\begin{equation*}
\m\beta(z)=\max\{0,\langle v,-n^*(z)\rangle\}, \m\m\m\m\m
n^*(z)=\argmax_{n\in n(z)}\langle v,-n\rangle.
\vspace{-2mm} 
\end{equation*}
\end{lemma}}

From Lemma~2.2, for any $z\in X$ and $v\in\mathbb{R}^q$, we have
\vspace{-2mm}
\begin{equation}\label{eq:Pi_omog}
\Pi_X(z,kv) = k\Pi_X(z,v) \FORALL k>0. \vspace{-5mm}
\end{equation} 

\begin{figure} \begin{center} 
   \includegraphics[height=3.25cm]{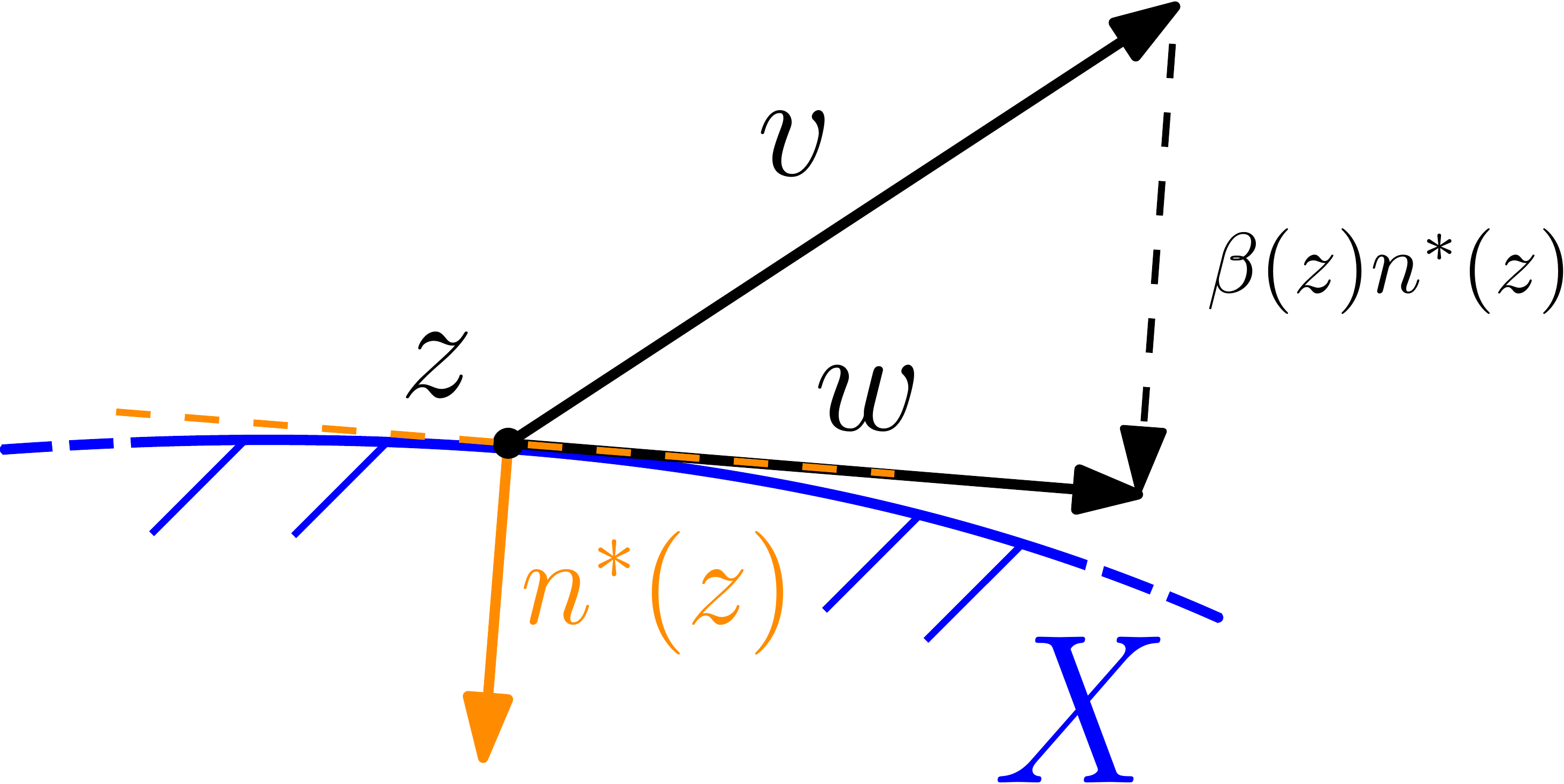}
   \caption{We show how $\Pi_X(z,\cdot)$ maps $v$ (pointing outward)
   to $w$, when $z\in\partial X$. (Here the set $n(z)$ is a singleton.)} \label{fig:proj} \end{center}
\end{figure}

To help the reader in understanding the intuition behind Lemma
\ref{lmm:Pi}, we show in Fig.~\ref{fig:proj} the resulting vector
$w=\Pi_X(z,v)$ when, e.g., $z\in\partial X$ and $v$ points outward.

\begin{definition} \label{def:sol}{\rm\cite[Definition 
   2.5]{Nagurney1995}.}
Let $X\subset\rline^q$ be a closed and convex set, $z\in X$, and $F:X
\to\rline^q$ a vector field. The function $z:[0,\infty)\to X$ is a
\textit{Carath\'{e}odory solution} to the equation \vspace{-1mm}
\begin{equation} \label{eq:PDS_ODE}
   \dot{z} \m=\m \Pi_X(z,-F(z)) \vspace{-1mm}
\end{equation} 
if $z(\cdot)$ is absolutely continuous and if \vspace{-2mm}
$$\dot{z}(t) \m=\m \Pi_X(z(t),-F(z(t))), \vspace{-2mm}$$ 
save on a set of Lebesgue measure zero (of points $t>0$).
\end{definition}

\vspace{-1mm}
For any $z_0\in X$ as initial value, we associate with
\rfb{eq:PDS_ODE} an initial value problem defined as: \vspace{-1mm}
\begin{equation} \label{eq:PDS_IVP}
   \dot{z} \m=\m \Pi_X(z,-F(z)), \qquad z(0)=z_0.
\end{equation}

\begin{remark}\label{rmk:in_X}
If \rfb{eq:PDS_IVP} has a solution, then such a solution is 
constrained in $X$ for all $t\geq0$.
\end{remark}

\begin{definition} \label{def:PDS}{\rm\cite[Definition 
   2.6]{Nagurney1995}.} 
Let $X$ and $F$ be as above. Define a \textit{projected dynamical
system} PDS($F,X$) as a map \m $\Phi:X\times\rline\mapsto X$, such
that $\phi_{z_0}(t)=\Phi(z_0,t)$ is a Carath\'{e}odory solution of
\rfb{eq:PDS_IVP}, so that for almost every $t>0$ \vspace{-1mm}
\begin{equation*}
   \dot{\phi}_{z_0}(t) \m=\m \Pi_X(\phi_{z_0}(t),-F(\phi_{z_0}(t))),
   \quad \phi_{z_0}(0)=z_0. \vspace{-4mm}
\end{equation*} 
\end{definition}

\begin{figure} \centering 
   \includegraphics[height=3.25cm]{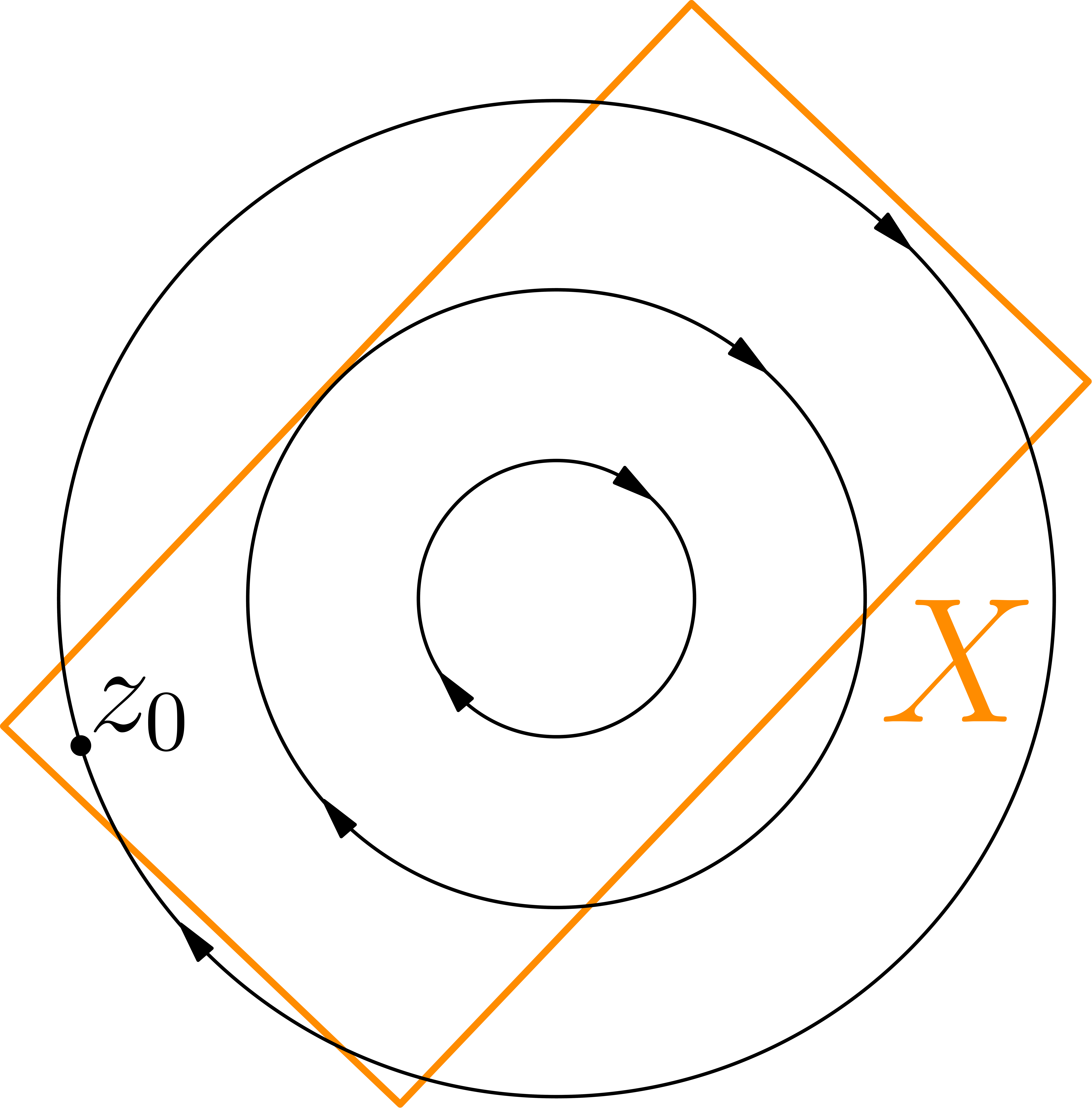}\qquad
   \includegraphics[height=3.25cm]{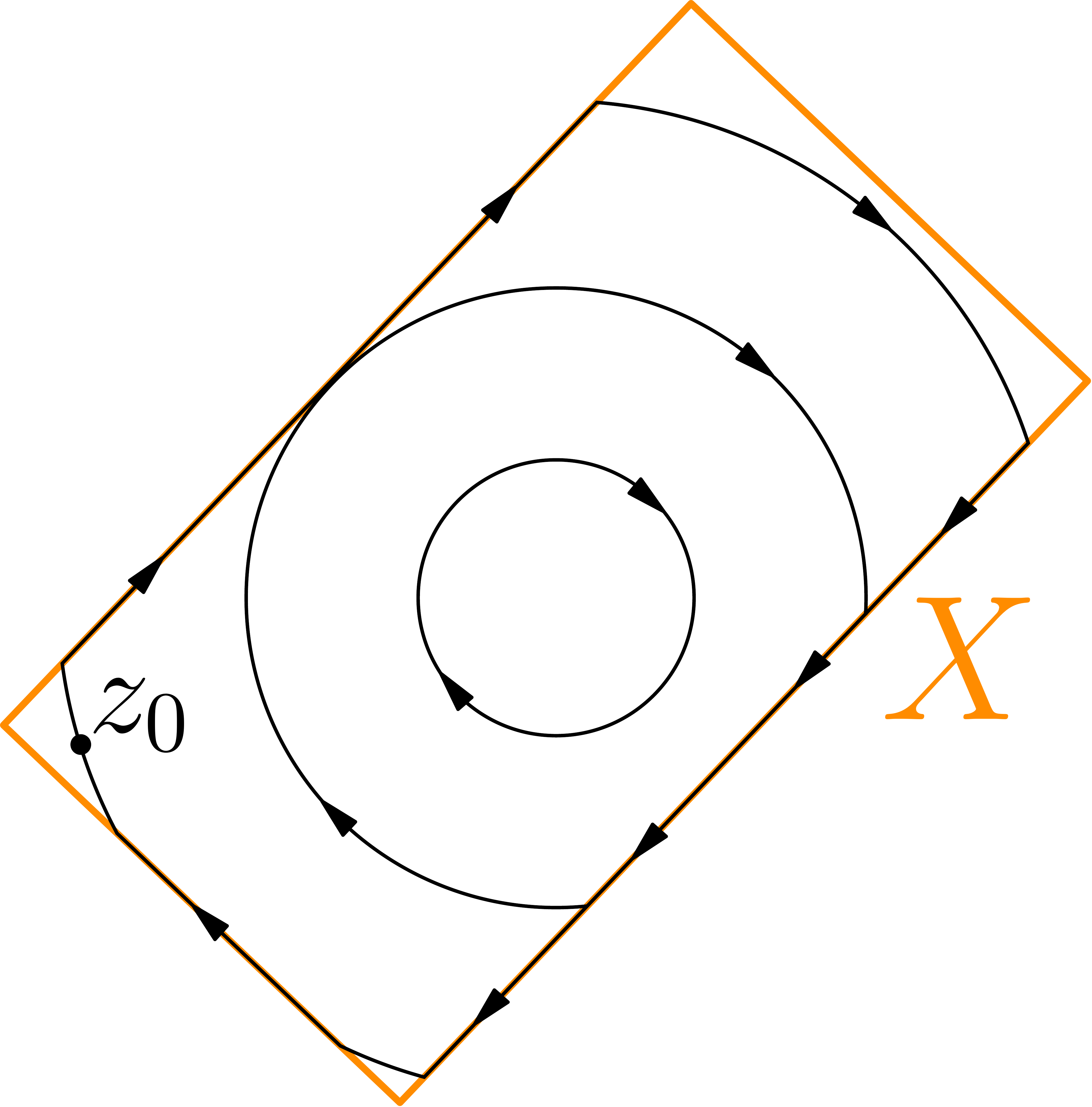} 
   \caption{Portrait of a classical dynamical system (left), and of a
   PDS corresponding to the same vector field (right). Adapted from
   \cite[Example 2.1]{Nagurney1995}.} \label{fig:PDS}
\end{figure}

We show in
Fig. \ref{fig:PDS} the portrait of a classical
dynamical system and the portrait of the corresponding PDS.

\begin{definition} \label{def:PDS_eq}{\rm\cite[Definition 
   2.7]{Nagurney1995}.} 
The vector $z^*\in X$ is an equilibrium point of the PDS($F,X$) if \m
$\Pi_X(z^*,-F(z^*))=0$.
\end{definition}

\begin{remark}
As pointed out in \cite{Nagurney1995}, $z^*\in X$ is
an equilibrium point of the PDS($F,X$) if the vector field $F$ 
vanishes at $z^*$. The converse, however, is only true when 
$z^*\in\mathrm{int}\m X$. Indeed, when 
$z^*\in\partial X$, we may have $F(z^*)\neq0$, but $\Pi_X(z^*,
-F(z^*))=0$.
\end{remark}

{\color{blue}\begin{theorem} \label{thm:ex_uniq}{\rm\cite[Theorem 
   2.5]{Nagurney1995}.} 
Assume that there exists a $B>0$ such that the vector field
$F:X\to\rline^q$ satisfies: \vspace{-1mm}
$$ \begin{gathered} \|F(z)\|\leq B(1+\|z\|) \ \ \forall  \ z\in X,
   \\ \langle -F(x) + F(y),x-y\rangle\leq B\|x-y\|^2 \ \ \forall \ 
   x,y\in X.\end{gathered}\vspace{-3mm}$$
Then: \vspace{-2mm}
\begin{enumerate}
\item For any $z_0\in X$, there exists a unique solution $z:[0,
   \infty) \to X$ to the initial value problem \rfb{eq:PDS_IVP}.
\item If $z_n\to z_0$ as $n\to\infty$, then $z(t;z_n)$ converges to 
   $z(t;z_0)$ uniformly on every compact set in $[0,\infty)$.
\end{enumerate}
\end{theorem}}

\begin{remark}\label{rmk:Pi_ext}
The definition of $\Pi_X$ can be extended for $z,v\in\rline^q$ as
follows: \vspace{-2mm}
\begin{equation} \label{eq:Pi_ext}
   \Pi_X(z,v) \m=\m \frac{P_X(z)-z}{\norm{P_X(z)-z}} \ \ \ \forall \
   z\in\rline^q\setminus X. \vspace{-1mm} 
\end{equation} 
Suppose that \rfb{eq:PDS_IVP} has a solution for any $z_0\in X$, given by
$\Phi$ from Definition \ref{def:PDS}. If $\tilde{z}_0\in\rline^q\setminus X$,
then the solution of \rfb{eq:PDS_IVP} (with $\Pi$ extended as in
\rfb{eq:Pi_ext} and $z(0)=\tilde{z}_0$)
will move with unit velocity towards $P_X(\tilde{z}_0)$,
until it reaches it (in finite time). Then it will follow the flow
$\Phi$. With this extension, Theorem \ref{thm:ex_uniq} remains valid
for all $z_0\in\rline^q$. \vspace{-1mm}
\end{remark}

\begin{remark} \label{rmk:Lip}
The (uniform) Lipschitz continuity of $F$ on $X\subset\rline^q$
implies the assumptions of Theorem \ref{thm:ex_uniq}. \vspace{-1mm}
\end{remark}


\begin{remark} \label{rmk:PDS_Hauswirth}
For the setting of this paper (i.e, $X$ closed and convex, and $F\in
C^1$), the theory on PDS developed in \cite{Nagurney1995} is
sufficient to derive our main result. However, for the interested
readers, we refer to the contribution
\cite{Hauswirth2021}, where the work of \cite{Nagurney1995}
is generalized in several directions. In particular, the conditions of
Theorem \ref{thm:ex_uniq} are relaxed and the existence and uniqueness
of Krasovskii (and, when possible, Carath\'{e}odory) solutions to
\rfb{eq:PDS_IVP} is proved under milder assumptions on the set $X$, on
the vector field $F$, and for a more general Riemannian metric, see
\cite[Table 1]{Hauswirth2021}. Using the equivalent formulation of
$\Pi_X$ from \cite{Hauswirth2021}, it can be checked
that $\Pi_X(z,\cdot)$ is a contraction. Indeed, defining the
\textit{tangent cone} $T_zX$ at $z\in X\subset\mathbb{R}^q$ as in
\cite[Definition~2.1]{Hauswirth2021}, with $X$ closed and convex, and
using as metric $g$ the Euclidean norm, then $\Pi_X(z,v)$ from
\rfb{eq:Pi} can be formulated as in
\cite[Definition~3.1]{Hauswirth2021}, i.e., \vspace{-1mm}
\begin{equation} \label{eq:Pi_Hauswirth}
   \Pi_X(z,v) \m=\m \argmin_{w\in T_zX} \|w-v\| \m. \vspace{-1mm}
\end{equation}
It is a well-known result that the above operator and, equivalently,
our \rfb{eq:Pi}, is a contraction, i.e., \vspace{-2mm}
\begin{equation} \label{eq:contraction}
   \norm{\Pi_X(z,v_1)-\Pi_X(z,v_2)}\leq\norm{v_1-v_2}, \vspace{-1mm}
\end{equation}
for all $v_1,v_2\in\rline^q$ and for all $z\in X$.
\end{remark}

\section{Problem formulation} \label{sec3} 

Consider the nonlinear plant $\mathbf{P_0}$ described by \vspace{-1mm}
\begin{equation} \label{eq:P}
   \dot{x} \m=\m f_0(x,v), \qquad y \m=\m g(x), \vspace{-1mm}
\end{equation} 
with $f_0\in C^2(\rline^n\times\Vscr;\rline^n)$, $g\in C^1(\rline^n;
\rline^p)$, where $\Vscr\subset\rline^m$ is an open
domain with $m\geq p$. \vspace{-1mm}

The control objective is to make the plant output signal $y$ track a
constant reference signal $r\in Y\subset\rline^p$, while making sure
that the plant input signal $v$ converges to a steady-state value in a
desired compact set $V\subset\rline^m$ (e.g., determined by
operational constraints). This tracking property should hold for all
plant initial states $x_0$ in a reasonably large open set in $\rline^n$.
\vspace{-1mm}

\begin{figure} \begin{center} 
   \includegraphics[width=0.40\textwidth]{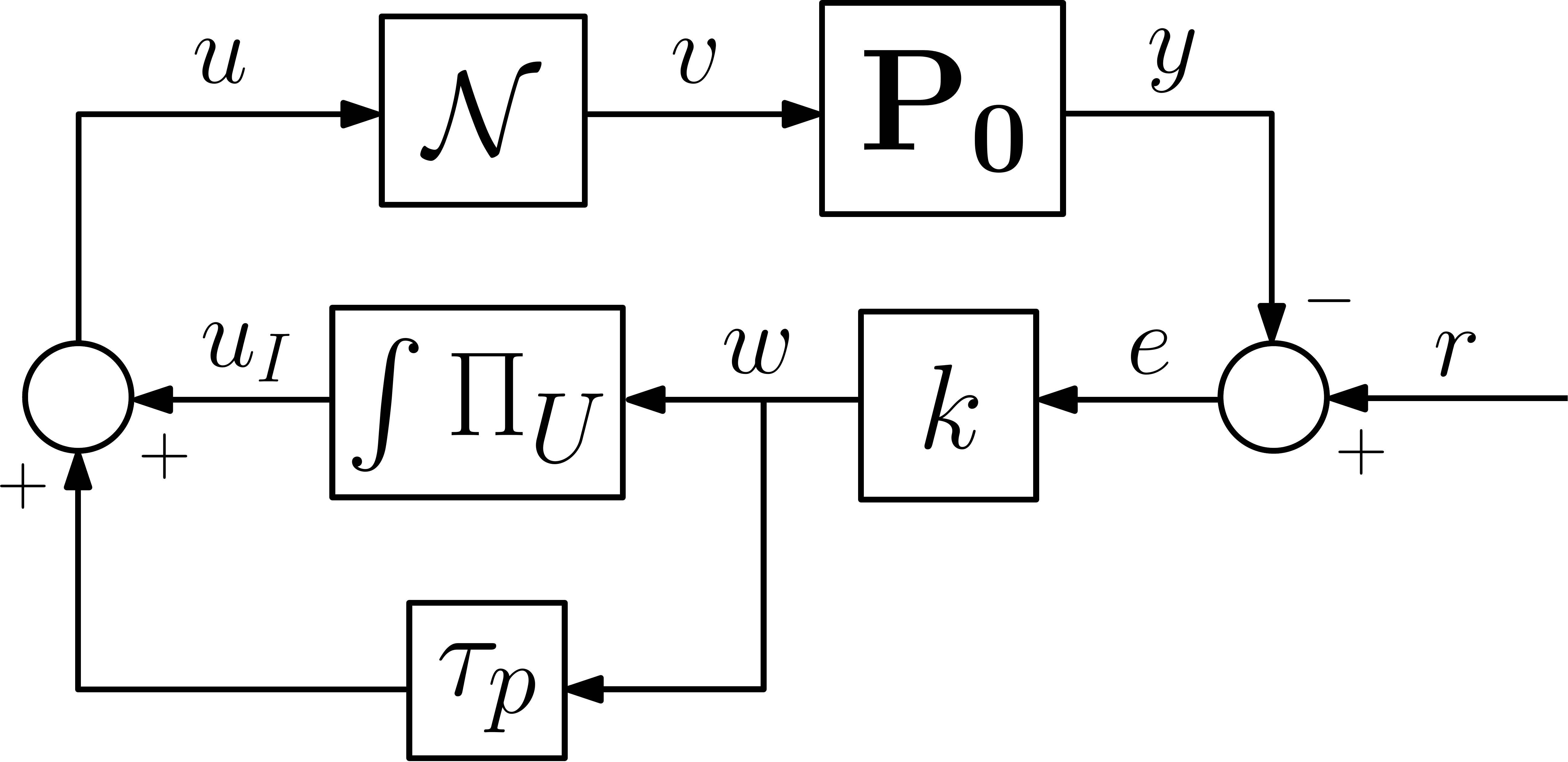}
   \caption{Representation of the closed-loop system 
   \rfb{eq:cl_MIMO_eq}, where $\mathbf{P_0}$ is the nonlinear plant 
   from \rfb{eq:P}.} \label{fig:clos_MIMO} \end{center}
\end{figure}

To achieve this control objective, we form the closed-loop system
shown in Fig.~\ref{fig:clos_MIMO}, described by the equations 
\vspace{-1mm}
\begin{equation} \label{eq:cl_MIMO_eq}
   \begin{gathered}
   \dot{x} \m=\m f_0(x,\Nscr(u_I+\tau_p k(r-g(x)))), \\
   \dot{u}_I \m=\m \Pi_U(u_I,k(r-g(x))), \vspace{-1mm}
   \end{gathered}
\end{equation}
where $\Uscr\subset\rline^p$ is an open domain, $U\subset\Uscr$ is a
compact and convex set ($\Uscr$ and $U$ to be defined), $\Pi_{U}$ is
the operator from \rfb{eq:Pi}, with the extension \rfb{eq:Pi_ext},
$\Nscr\in C^2(\Uscr,\Vscr)$ (to be defined), $V=\Nscr(U)$, $k>0$ and
$\tau_p\geq 0$. The {\em state space} of \rfb{eq:cl_MIMO_eq} is 
$\rline^n\times\Uscr$ and its {\em state} is $z(t)= \sbm{x(t)\\ 
u_I(t)}$. As is often the case in nonlinear systems theory, the first 
equation in \rfb{eq:cl_MIMO_eq} only makes sense on a ``region of 
interest'' in the state space, namely on the open set \vspace{-2mm}
$$ \Dscr_r \m\vcentcolon=\m \left\{ \sbm{x\\ u_I}\in \rline^n\times 
   \Uscr \ \big|\ u_I+\tau_p k(r-g(x)) \in \Uscr \right\} \m. 
   \vspace{-5mm}$$ 

It will be convenient to introduce the ``new plant'' $\PPP$ as the 
cascade of \m $\Nscr$ and $\mathbf{P_0}$, described by \vspace{-2mm}
\begin{equation} \label{eq:f}
   \dot{x} \m=\m f(x,u), \qquad y \m=\m g(x), \vspace{-1mm}
\end{equation} 
where $f(x,u)\vcentcolon=f_0(x,\Nscr(u))\in C^2(\rline^n\times
\Uscr;\rline^n)$.


{\color{blue}
\begin{proposition} \label{prop:uniqueness}
Consider the closed-loop system \rfb{eq:cl_MIMO_eq}, with $k,\tau_p
\in\rline$, $r\in\rline^p$. Then for every $\sbm{x_0\\ u_0}\in\Dscr_r$
with $u_0\in U$, there exists $\tau\in(0,\infty]$ such that
\rfb{eq:cl_MIMO_eq}, with initial conditions $z(0)=\sbm{x_0\\ u_0}$,
has a unique Carath\'{e}odory solution (or state trajectory) $z=\sbm
{x\\ u_I}$ defined on $[0,\tau)$. If $\tau$ is finite and maximal 
(i.e., the state trajectory cannot be continued beyond \m $\tau$),
then $\limsup_{t\rarrow\tau}\|x(t)\|=\infty$, or the signal $u(t)=
u_I(t)+\tau_p k(r-g(x(t)))$ approaches $\partial\Uscr$: \vspace{-2mm}
\begin{equation} \label{Ben_and_Jerry}
   \liminf_{t\rarrow\tau}d(u(t),\partial\Uscr) \m=\m 0 \m, 
   \vspace{-2mm}
\end{equation}
where $d$ denotes the distance in $\rline^p$. \vspace{-3mm}
\end{proposition}}

\begin{pf}
We introduce the closed and convex set $X\vcentcolon=\rline^n
\times U$. An equivalent representation of \rfb{eq:cl_MIMO_eq} (for 
$u_I(t)\in U$) is \vspace{-2mm}
\begin{equation} \label{eq:cl_MIMO}
   \dot{z} \m=\m \Pi_X(z,-F(z)), \vspace{-1mm}
\end{equation}
where \vspace{-1mm}
$$ -F(z) \m\vcentcolon=\m \sbm{f(x,u_I+\tau_p k(r-g(x))) \\ 
   k(r-g(x))},$$
so that \rfb{eq:cl_MIMO} makes sense as long as $z(t)\in\Dscr_r$. For
any $\delta>0$, $B_\delta$ denotes the closed ball of radius $\delta$
in $\rline^n$, and also in $\rline^p$ (the dimension will be clear 
from the context). We fix $\sbm{x_0\\ u_0}\in\Dscr_r$ such that 
$u_0\in U$. Define \vspace{-2mm}
$$ X_\delta \m\vcentcolon=\m (x_0+B_\delta)\times \left[ (u_0
   +B_\delta)\cap U \right] \m. \vspace{-2mm}$$
We choose $\delta$ small enough so that $X_\delta\subset\Dscr_r$. We
have $F\in C^1(\Dscr_r,\rline^n\times\rline^p)$ and $X_\delta$ is
compact and convex, thus it follows from Theorem \ref{thm:ex_uniq}
and Remark \ref{rmk:Lip} that \rfb{eq:cl_MIMO}, but with $X_\delta$ 
in place of $X$, has a unique solution $z:[0,\infty)\to X_\delta$
that satisfies $z(0)=\sbm{x_0\\ u_0}$. As long as $\|z(t)-z(0)\|<
\delta$, this solution $z$ is also a solution of the original
\rfb{eq:cl_MIMO}. From here it follows that there exists $\tau>0$
such that \rfb{eq:cl_MIMO} has a unique state trajectory $z$ defined
on $[0,\tau)$, starting from the initial state $z(0)=\sbm{x_0\\ 
u_0}$. \vspace{-2mm}

Suppose that $\tau>0$ as above is
finite and maximal. If $\limsup_{t\rarrow\tau}\|x(t)\|$ is finite,
then the trajectory $z$ is bounded on $[0,\tau)$ (because $u_I(t)\in
U$ for all $t\in[0,\tau)$). If \rfb{Ben_and_Jerry} were not true,
then there exists $\e>0$ such that $d(u(t),\partial\Uscr)\geq\e$ for
all $t\in[0,\tau)$. This implies that the closure of $\{z(t)\ |\ t
\in[0,\tau)\}$ is a compact subset of $\Dscr_r$. Since $F$ is 
continuous on $\Dscr_r$, there exists $M>0$ such that $\|F(z(t))\|
\leq M$ for all $t\in[0,\tau)$. Let $(t_j)$ be an increasing 
sequence such that $t_j\in[0,\tau)$, $t_j\rarrow\tau$. Using
\rfb{eq:contraction} we obtain that for $j>k$ \vspace{-3mm}
$$ \|z(t_j)-z(t_k)\| \m\leq\m \int_{t_k}^{t_j}\|\Pi_X(z(t),-F(z(t)))
   \| \dd t \vspace{-4mm}$$
$$ \m\hspace{24mm} \m\leq\m M(t_j-t_k) \m.$$
Thus, $(z(t_j))$ is a Cauchy sequence, so that it converges to a 
limit $z(\tau)\in\Dscr_r$. It is easy to see that the limit is 
independent of the choice of $(t_j)$, and 
that the function $z$, extended to $[0,\tau]$, is a 
Carath\'{e}odory solution of \rfb{eq:cl_MIMO} on $[0,\tau]$. We
could extend this solution even further, using the argument in the
first part of this proof. This would contradict the maximality of 
$\tau$, hence our assumption that \rfb{Ben_and_Jerry} is false has
led us to a contradiction. Thus, if $\limsup_{t\rarrow\tau}\|x(t)
\|$ is finite, then \rfb{Ben_and_Jerry} holds. \hfill $\blacksquare$     
\end{pf}

\vspace{-5mm}
\begin{remark}\label{rmk:Haws}
The closed-loop system \rfb{eq:cl_MIMO_eq}, with $\tau_p=0$ and
$\Nscr=I$, can be approximated by the equations \vspace{-2mm}
\begin{equation} \label{eq:cl_Hauswirth}
   \begin{gathered} \dot{x}=f(x,P_U(u_I)), \\ \dot{u}_I=k(r-g(x))-
   \frac{1}{K}(u_I-P_U(u_I)), \vspace{-2mm} \end{gathered} 
\end{equation}
with $K>0$ small and $P_U$ from \rfb{eq:proj}. Indeed, for every initial state 
$(x_0,u_0)\in \mathbb{R}^n\times U$, the solution of \rfb{eq:cl_Hauswirth} converges
uniformly to that of \rfb{eq:cl_MIMO_eq} (with $\tau_p=0$, and 
$\Nscr=I$) for $K\to0^+$, see \cite[Theorem~2]{Hauswirth2020b}.
\end{remark}

\section{Closed-loop stability analysis} \label{sec4} 

In this section we present our main result, namely, we derive an upper
bound for the gain $k$ ensuring the existence of a (locally)
exponentially stable equilibrium point for the closed-loop system
\rfb{eq:cl_MIMO_eq}, for each constant reference 
$r\in Y\subset\mathbb{R}^p$ ($Y$ to be defined). 
We further characterize a subset
of the region of attraction of this equilibrium point 
such that if the initial state is in this region, then
the plant output $y$ tracks $r$. 
This result generalizes \cite[Theorem~4.3]{Lor2020}, which
was formulated for the SISO saturating integrator. As in
\cite{Lor2020}, our stability analysis employs SP
methods (see Appendix \ref{app:A} for the details), which
can be found, e.g., in \cite[Ch.~7]{Kokotovic1999}, \cite[Ch.~11]
{Khalil2002}.

\vspace{2mm}\begin{mdframed}
\begin{assum}\label{assumption_1}
There exists a function \\  $\Xi\in C^1(\Vscr;\rline^n)$ such that
\vspace{-2mm}
\begin{equation} \label{eq:xi_equil}
   f_0(\Xi(v),v) \m=\m 0 \FORALL v \in \Vscr. \vspace{-2mm}
\end{equation}
Moreover, the equilibrium points $\{\Xi(v)\m\big|\m v\in
\Vscr\}$ are uniformly exponentially stable. This means that there
exist $\e_0>0$, $\l>0$ and $\rho\geq 1$ such that for each constant
input $v_0\in\Vscr$, the following holds: \vspace{-3mm}

If $\norm{x(0)-\Xi(v_0)}\leq\e_0$, then for every $t\geq 0$,
\vspace{-2mm}
\begin{equation} \label{eq:exp_stab}
   \norm{x(t)-\Xi(v_0)} \leq \rho e^{-\l t}\norm{x(0) - \Xi(v_0)}.
\end{equation}
\end{assum}
\end{mdframed}

\vspace{-1mm}
\begin{remark}
Assumption \ref{assumption_1} guarantees the stability of the
boundary-layer system associated to the closed-loop system
\rfb{eq:cl_MIMO_eq} (see \rfb{eq:fast_system} in Appendix \ref{app:A}). This is a standard
assumption in the framework of SP theory (see, for
instance, \cite{Desoer1985}, \cite[Ch.~11]{Khalil2002}, \cite
[Ch.~7]{Kokotovic1999}).
%
\vspace{-1mm}
\end{remark}

\begin{remark} \label{rmk:jacobian}
The (uniform) exponential stability condition \rfb{eq:exp_stab} can be
checked by linearization: If the Jacobian matrices \vspace{-2mm} 
$$ A(v_0) \m=\m \left.  \frac{\partial f_0(x,v)}{\partial x} \right|_
   {\nm\begin{array}{c} \m \scriptstyle x\m=\m\Xi(v_0)\vspace{-3mm}\\
   \scriptstyle v\m=v_0 \end{array}} \m\in\m \rline^{n\times n}
   \vspace{-5mm}$$
have eigenvalues bounded away from the right half-plane, \vspace{-2mm}
$$ \max\Re\sigma(A(v_0)) \m\leq\m \l_0 \m<\m 0 \FORALL v_0\in\Vscr,
   \vspace{-2mm}$$
then $\Xi(v_0)$ is a uniformly exponentially stable equilibrium point of $\m\mathbf{P_0}$,
for all $v_0\in\mathcal{V}$, see
\cite[eq.~(11.16)]{Khalil2002}.
\end{remark}

\vspace{-1mm}
\begin{remark}\label{rmk:ex_Xi}
If $\Xi$ satisfies \rfb{eq:xi_equil}
and \rfb{eq:exp_stab}, then $\Xi\in C^2$ thanks to
the implicit function theorem (since $f_0\in C^2$).
\end{remark}

\vspace{-1mm}
\textbf{Notation.} Let $G(v)\vcentcolon=g(\Xi(v))\in C^1(\Vscr;
\rline^p)$ denote the steady-state input-output map corresponding to
$\mathbf{P_0}$.

\vspace{2mm}\begin{mdframed}
\begin{assum} \label{assumption_2}
The plant $\mathbf{P_0}$ satisfies Assumption \ref{assumption_1}.
Moreover, there exist an open domain $\Uscr\subset\rline^p$, a function
$\Nscr\in C^2(\Uscr,\Vscr)$, and $\mu>0$ such that \vspace{-2mm}
$$ \langle\m G(\Nscr(u_1))-G(\Nscr(u_2)),u_1-u_2\m\rangle \m\geq \mu\m
   \norm{u_1-u_2}^2 \vspace{-1mm}$$
for all $u_1,u_2\in\Uscr$, i.e., $G\circ\Nscr$ is strictly monotone.
\end{assum} \end{mdframed}

\vspace{-1mm} 
We choose $U\subset\Uscr$ to be compact, convex, with
$\mathrm{int}\m U\neq\emptyset$. We let $Y=G(\Nscr(U))$, and, for 
any $r\in Y$, we define \vspace{-2mm}
$$ u_r\vcentcolon \m=\m (G\circ\Nscr)^{-1}(r) \qquad x_r\vcentcolon
   \m=\m \Xi(\Nscr(u_r)),\vspace{-2mm}$$
which are well-defined since $G\circ\Nscr$ is strictly monotone on
$\mathcal{U}$ (hence one-to-one). From Assumption~\ref{assumption_1},
$(x_r,u_r)$ is an equilibrium point
of the closed-loop system \rfb{eq:cl_MIMO_eq}.

\vspace{-2mm} \textit{Some commentary on the sets $\Vscr$, $\Uscr$,
$V$, $U$, and $Y$.} The set $\Vscr\subset\rline^m$ is a set of inputs
for which we have steady-state stability of the plant $\mathbf{P_0}$
(see Assumption \ref{assumption_1}). The set $\Uscr\subset\rline^p$ is
where Assumption \ref{assumption_2} holds and, thus, where we would
like to constrain the state of the integrator $u_I$ in order to obtain
closed-loop stability. The set $\Uscr$ may be too large, and, to
satisfy operational constraints, we impose $u_I(t)\in U$, where $U$ is
chosen as above. We denote $V=\Nscr(U)$. Finally, $Y=G(\Nscr(U))$ is
the natural set of feasible references, since $y=G(\Nscr(u_I))$ at 
steady-state.

\vspace{-2mm}
\begin{remark}
Assumption~\ref{assumption_2} guarantees the stability of the reduced-order
model associated to the closed-loop system \rfb{eq:cl_MIMO_eq} (see \rfb{eq:reduced_model} in Appendix
\ref{app:A} for the details). This is a common assumption
when SP tools are used to investigate the stability
of a nonlinear plant connected in feedback with an
integral controller, see \cite{Desoer1985,Huang2019}. The work \cite{Simpson-Porco2020} has extended
the result from \cite{Desoer1985}, by replacing the monotonicity
assumption on the input-output steady-state map with the 
infinitesimal contracting property of the reduced dynamics.
However, as discussed in \cite[Sect.~3]{Simpson-Porco2020}, if the 
infinitesimal contracting property is stated with respect to the standard Euclidean norm, then
the conditions of \cite{Simpson-Porco2020} reduce to those of \cite{Desoer1985}. 
In our framework, the two are equivalent.
\vspace{-1mm}
\end{remark}

\vspace{-1mm}
\begin{remark}\label{rmk:G_inc}
For a matrix $M\in\rline^{p\times p}$, define $\Re M=\half (M+M^\top)
$. The (strict) monotonicity of $G\circ\Nscr\in C^1(\Uscr,\rline^p)$
is equivalent to the fact that $\Re\frac{\partial(G\circ\Nscr)}
{\partial u}$ is strongly positive, i.e., there exists a $\mu>0$ such
that \vspace{-1mm}
$$ \left\langle\frac{\partial (G\circ\Nscr)}{\partial u}w,w\right
   \rangle\geq\mu\norm{w}^2 \ \ \ \ \forall \ w\in\rline^p, \ \forall
   \ u\in \Uscr,\vspace{-1mm}$$
see \cite[Proposition~2.5]{Nagurney1995}.
\end{remark} 

\vspace{-2mm}
There are several ways to choose $\Nscr$ (the case $\Nscr=I$ was
considered in our recent conference paper \cite{Lor2021Conf}). 
Assume that $G$ admits a right inverse $G^{-1}_{\mathrm{right}}
\in C^2(G(\Vscr);\Vscr)$, i.e., $G\circ G^{-1}_{\mathrm{right}}=I$ 
(the identity on $\Vscr$). 
Then, we suggest the choice $\Nscr= G^{-1}_{\mathrm{right}}$,
for which Assumption \ref{assumption_2} trivially holds, and
$U=Y$. 
As shown in Sect.~\ref{sec5}, the choice 
$\Nscr=G^{-1}_{\rm right}$ can be very convenient.

\vspace{-1mm}
\begin{remark}
If $\mathbf{P_0}$ is linear, described by the matrices $A,B,C$ in the
usual way ($\dot{x}=Ax+Bv$, $y=Cx$), then Assumption
\ref{assumption_1} reduces to the fact that $A$ is Hurwitz. The
functions $\Xi,G$ from Assumptions \ref{assumption_1} and
\ref{assumption_2} are given by \vspace{-2mm}
$$ \Xi(v) \m=\m (-A)^{-1}Bv, \qquad G(v) \m=\m P(0)v, \vspace{-2mm}$$
where $P(s)=C(sI-A)^{-1}B$ is the plant transfer function. In this
case, if $P(0)$ is onto, then $\Nscr$ can be chosen as
$\Nscr=P(0)^*(P(0)P(0)^*)^{-1}$, so that, again, $G\circ\Nscr=I$, and
Assumption \ref{assumption_2} is trivially satisfied.
\end{remark}


\vspace{-2mm}
{\color{blue}
\begin{theorem} \label{thm:cl_stab}
Consider the closed-loop system \rfb{eq:cl_MIMO_eq}, where
$\mathbf{P_0}$ satisfies Assumption \ref{assumption_2}. Then there
exists a $\kappa>0$ such that if the gain $k\in(0,\kappa]$, then for
any $r\in Y=G(\Nscr(U))$, $(\Xi(\Nscr(u_r)),u_r)$ is a (locally)
exponentially stable equilibrium point of the closed-loop system
\rfb{eq:cl_MIMO_eq}, with state space $\mathcal{X}=\rline^n\times
\mathcal{U}$. If the initial state $\sbm{x_0\\ u_0}\in \mathcal{D}_r$
($\mathcal{D}_r$ from Sec.~3)
of the closed-loop system satisfies $u_0\in U$ and 
$\|x_0-\Xi(\mm\Nscr(u_0))\|\leq\e_0$, then \vspace{-3mm}
\begin{equation}\label{eq:conv}
x(t)\rarrow\Xi(\Nscr(u_r)), \qquad u_I(t)\rarrow u_r, \qquad y(t)
\rarrow r, \vspace{-3mm}
\end{equation}
and this convergence is at an exponential rate. \vspace{-2mm}
\end{theorem}}

For the proof see Appendix \ref{app:A}. Note that, clearly,
\rfb{eq:conv} implies that $u(t)\rarrow u_r$ (since $e(t)=r-g(x(t))\rarrow0$).
\vspace{-2mm}
%

\begin{remark} \label{rmk:global}
The results from Theorem~\ref{thm:cl_stab} can be extended
globally, following the procedure of \cite[Sect.V]{Lor2020},
if $\mathbf{P_0}$ satisfies the asymptotic gain
property (introduced in \cite{Sontag1996}) around each equilibrium
point $\Xi(v_0)$, for all $v_0\in\Vscr$.
\end{remark}

\section{Power regulation for a grid-connected synchronverter}
         \label{sec5} 

We present an application of the proposed control strategy for the
(active and reactive) power regulation of a grid-connected
synchronverter, when the grid is modelled as an infinite bus. In our
simulations, we assume that the power set points for the
synchronverter control algorithm are provided by an external control
loop (e.g., using optimal power flow considerations), which we do not model. The synchronverter output active and
reactive powers have to track these set points, whilst making sure to
not leave the safe operating region. We compare the behaviour of the
closed-loop system formed by the synchronverter model $\mathbf{P_0}$,
our saturating integrator $\int\Pi_{U}$ (here $\tau_p=0$), and the
nonlinear gain $\Nscr=G^{-1}_{\mathrm{right}}$ (to be defined), 
with the one formed by $\mathbf{P_0}$, a classical integrator 
($\Pi_{U}=I,\m\tau_p=0$), and a static linear gain $\Nscr=K\in
\rline^{p\times p}$ (to be defined). \vspace{-1mm}

\subsection{Description of the synchronverter model} \vspace{-1mm}

Synchronverters, see \cite{Zhong2011}, are a particular type
of \textit{virtual synchronous machines}, i.e., inverters with a
control algorithm that causes them to behave towards the power grid
like synchronous generators.  Among the different grid-connected
synchronverter models in the literature, we refer to the fourth order grid-connected
synchronverter model from \cite[eq.~(3.1)]{Natarajan2017},
\cite[eq.~(13)]{Lor2021}, where the grid is modelled as an infinite
bus.  Due to lack of space, we omit the physical meaning of the
equations, which can be found in the just cited references.

\vspace{-1mm}
Let $\mathbf{P_0}$ be the fourth order grid-connected synchronverter
model with state \vspace{-2mm}
\begin{equation} \label{eq:x} 
   x \m=\m [ \ i_d \ i_q \ \o \ \delta \ ]^\top \m \in \m \rline^4,
   \vspace{-2mm}
\end{equation} 
where $i_d$ and $i_q$ are the d and q components of the stator
currents, $\o$ is the (virtual) rotor angular velocity, and $\delta$
is the power angle (regarded modulo $2\pi$, i.e., $\delta$ and
$\delta+2 \pi$ are considered to be the same angle). The input is
\vspace{-2mm}
\begin{equation} \label{eq:u} 
   v \m=\m [ \ T_m \ i_f \ ]^\top \m \in\m \rline\times(0,\infty), 
   \vspace{-2mm}
\end{equation} 
where $T_m$ is the (virtual) prime mover torque, and $i_f$ is the
(virtual) field current. The output is \vspace{-2mm} 
\begin{equation} \label{eq:y} 
   y \m=\m [ \ P \ Q \ ]^\top \m \in\m \rline^2, \vspace{-2mm}
\end{equation} 
where $P$ is the active power, and $Q$ is the reactive power. The
plant $\mathbf{P_0}$ is described by the equations \vspace{-2mm}
\begin{equation} \label{eq:ss4}
   \begin{gathered} H\dot{x} \m=\m A(x,v)x \m + h(x,v) \m, \\ 
   y = g(x), \end{gathered}\vspace{-2mm} 
\end{equation} 
with \vspace{-2mm} 
$$ H = \bbm{\ L & \ 0 & \ 0 & \ 0  \ \\ \ 0 & \ L & \ 0 & \ 0 \ \\
   \ 0 & \ 0 & \ J & \ 0 \ \\ \ 0 & \ 0 & \ 0 & \ 1 \ }\m,\quad \m 
   h(x,v) \m= \bbm{\ V\sin \delta \ \\ \ V\cos\delta \ \\ \ T_m \m + 
   D_p\o_n \ \\ \ \vspace{1mm} -\o_g \ }, \vspace{-1mm}$$
$$ A(x,v) \m=\m \bbm{\ -R & \ \o L & \ 0 & \ 0 \ \\ \ -\o L & \ -R & 
   \ -mi_f & \ 0 \ \\ \ 0 & \ m i_f & \ -D_p & \ 0 \ \\ \ 0 & \ 0 & 
   \ 1 & \ 0 \ },$$
and
$$ g(x) \m=\m -V\bbm{\ \cos\delta & \ \sin\delta \ \\ \ -\sin\delta
   & \ \cos\delta\ } \bbm{\ i_q \ \\ \ i_d \ } \m.$$
Here $L>0$ is the total stator inductance, $R>0$ is the total stator
resistance, $J>0$ is the rotor moment of inertia, $V>0$ is the rms 
value of the line voltage, $D_p>0$ is the frequency droop constant,
$\o_n$ is the nominal grid frequency, $\o_g$ is the grid frequency, 
and $m=\sqrt{3/2}M_f$, where $M_f>0$ is the peak mutual inductance
between the virtual rotor winding and any one stator winding.

\vspace{-1mm}
\textit{Synchronverter parameters.} We use the synchronverter
parameters from the numerical example \cite[Subsect.~VI-A]{Lor2021}, 
chosen for a synchronverter designed to supply a nominal active power
of $9$\m kW to a grid with nominal frequency $\o_n=100\pi$\m rad/sec (50\m Hz) 
and line voltage $V=230\sqrt{3}$\m Volts. The parameters are: $J=0.2
$\m Kg$\cdot$m$^2$/rad, $D_p=3$\m N$\cdot$m/(rad/sec), $R=1.875$\m 
$\Om$, $L= 56.75$\m mH, $m=3.5$\m H, and $\o_g=\o_n$. \vspace{-1mm}

\vspace{-1mm}
\subsection{Formulation of the control problem}\vspace{-1mm} 
The control problem that we address is the regulation of the synchronverter output
$y$ to the reference signal \vspace{-2mm}
\begin{equation} \label{eq:r} 
   r \m=\m [ \ P_{\rm set} \ Q_{\rm set} \ ]^\top \m\in\m \rline^2,
   \vspace{-2mm}
\end{equation}
while keeping the synchronverter input $v$ in a safe (compact and convex)
operating region $V\subset\rline^2$ ($m=p=2$). We form a
closed-loop system as in Fig.~\ref{fig:clos_MIMO} (here $\tau_p=0$),
and we are interested in studying its stability and tracking
properties using Theorem \ref{thm:cl_stab}. To this aim, we first
verify whether $\mathbf{P_0}$ from \rfb{eq:ss4} satisfies Assumptions~\ref{assumption_1},~\ref{assumption_2}.

\vspace{-1mm} \textit{Verification of Assumption \ref{assumption_1}.}
The equilibrium points of the grid-connected synchronverter model \rfb{eq:ss4}
have been studied in \cite{Lor2021}. In particular, in \cite[Prop.~3.1 and
Prop.~3.3]{Lor2021} it is shown that for each \m $T_m>-\frac{V^2}{4R\o_g}$ there
is a finite interval $I_f\subset(0,\infty)$ such that for $i_f\in{\rm int}\m I_f$,
the model $\mathbf{P_0}$ has two equilibrium points, of which at most one is
stable. 
We denote by $\Vscr$ the subset of $\rline\times(0,\infty)$ such that if 
$v\in\Vscr$, then $\mathbf{P_0}$ has an
exponentially stable equilibrium point corresponding to the constant input
signal $v$.

We mention that in \cite[Theorem 6.3]{Natarajan2018} sufficient 
conditions were given for an equilibrium point of $\mathbf{P_0}$ to
be almost globally asymptotically stable. In this paper, (local)
exponential stability is what we need, and that can be checked with 
relative ease, using the linearization of \rfb{eq:ss4}, according to
Remark \ref{rmk:jacobian}. \vspace{-1mm}

\begin{figure} \begin{center} 
   \includegraphics[height=5cm]{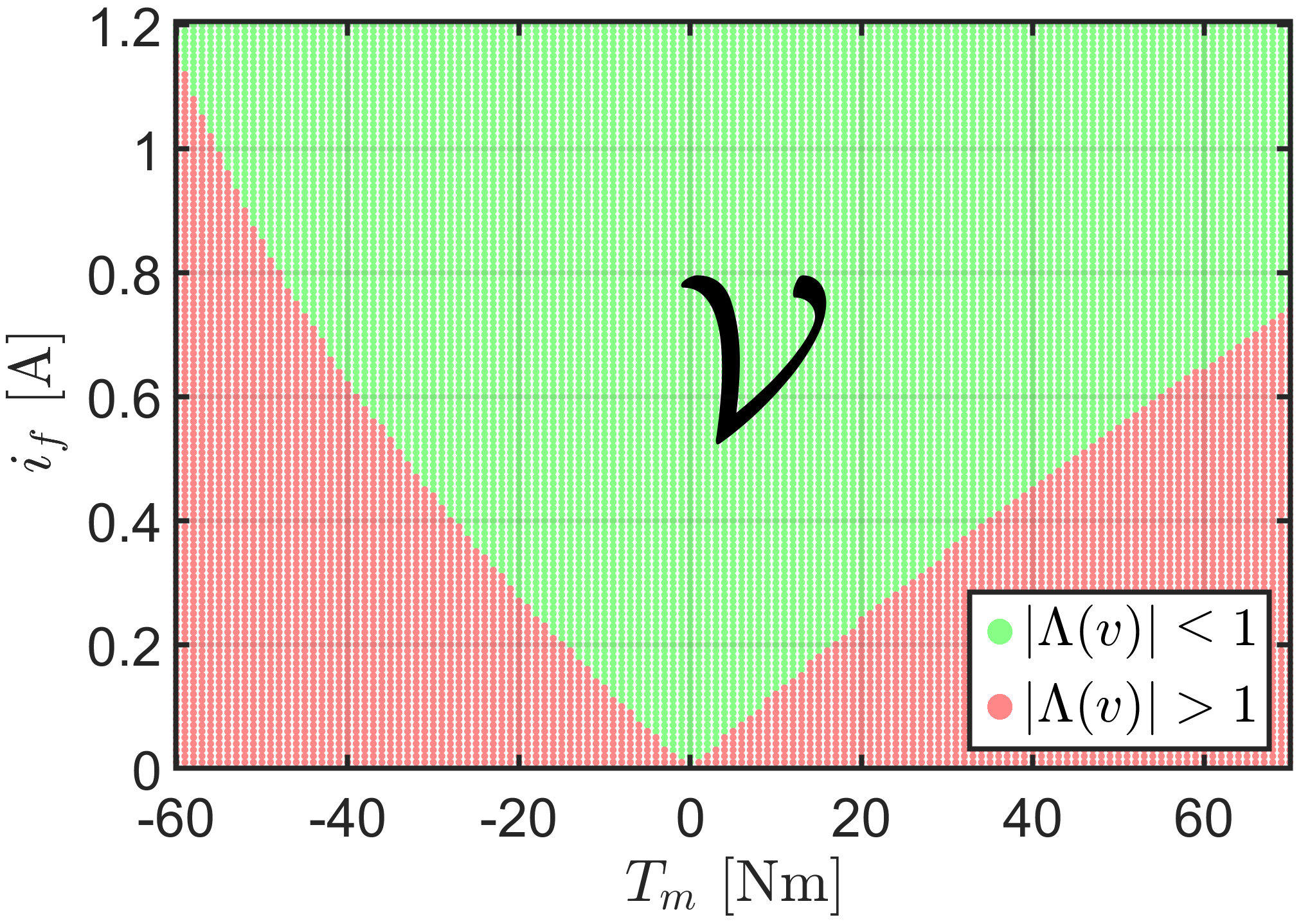}
   \caption{A glance of the set $\Vscr$ (in light green), when
   $[T_m \ i_f]^\top\in\mathcal{R}=[-60,70]\times[0.01,1.2]$.} \label{fig:V}
   \end{center}
\end{figure}

The function $\Xi:\Vscr\to\rline^4$ is given by \vspace{-1mm}
$$ \Xi(v) \m=\m \begin{bmatrix} \ -\frac{T_m\o_g}{m i_fp}+\frac{V
   \sin(\arccos\L(v)-\phi)}{R} \ \\ \ -\frac{T_m}{mi_f} \ \\ \ 
   \o_g \ \\ \ \arccos\L(v)-\phi \ \end{bmatrix},$$ 
where \m $\phi\in\left(0,\frac{\pi}{2}\right)$ such that \m $\tan
\phi=\frac{\o_gL}{R}$, \vspace{-3mm}
$$ \L(v) \m=\m -\frac{T_m}{mi_f} \frac{L\sqrt{p^2+\o_g^2}} {V} + 
   \frac{mi_f\o_g p}{V\sqrt{p^2 + \o_g^2}}, \ \ \ p=\frac{R}{L}.$$
Our numerical explorations indicate that within the rectangle $\Rscr:=[-60,70]
\times[0.01,1.2]$
shown in Fig.~\ref{fig:V}, 
\vspace{-2mm}
$$ \mathcal{V}\cap\mathcal{R}=\{\m[T_m \ i_f]^\top\in\Rscr \ \big|\ \m |\L(v)|<1\m\}. \vspace{-2mm}$$
For more details 
on the mapping $\Xi$, see \cite[Proposition~3.1 and eq.~(48)]{Lor2021}
($\tilde{T}_m$ there is equivalent to $T_m$ here, since we have assumed $\o_n=\o_g$).


\begin{figure} \begin{center} 
   \includegraphics[height=5cm]{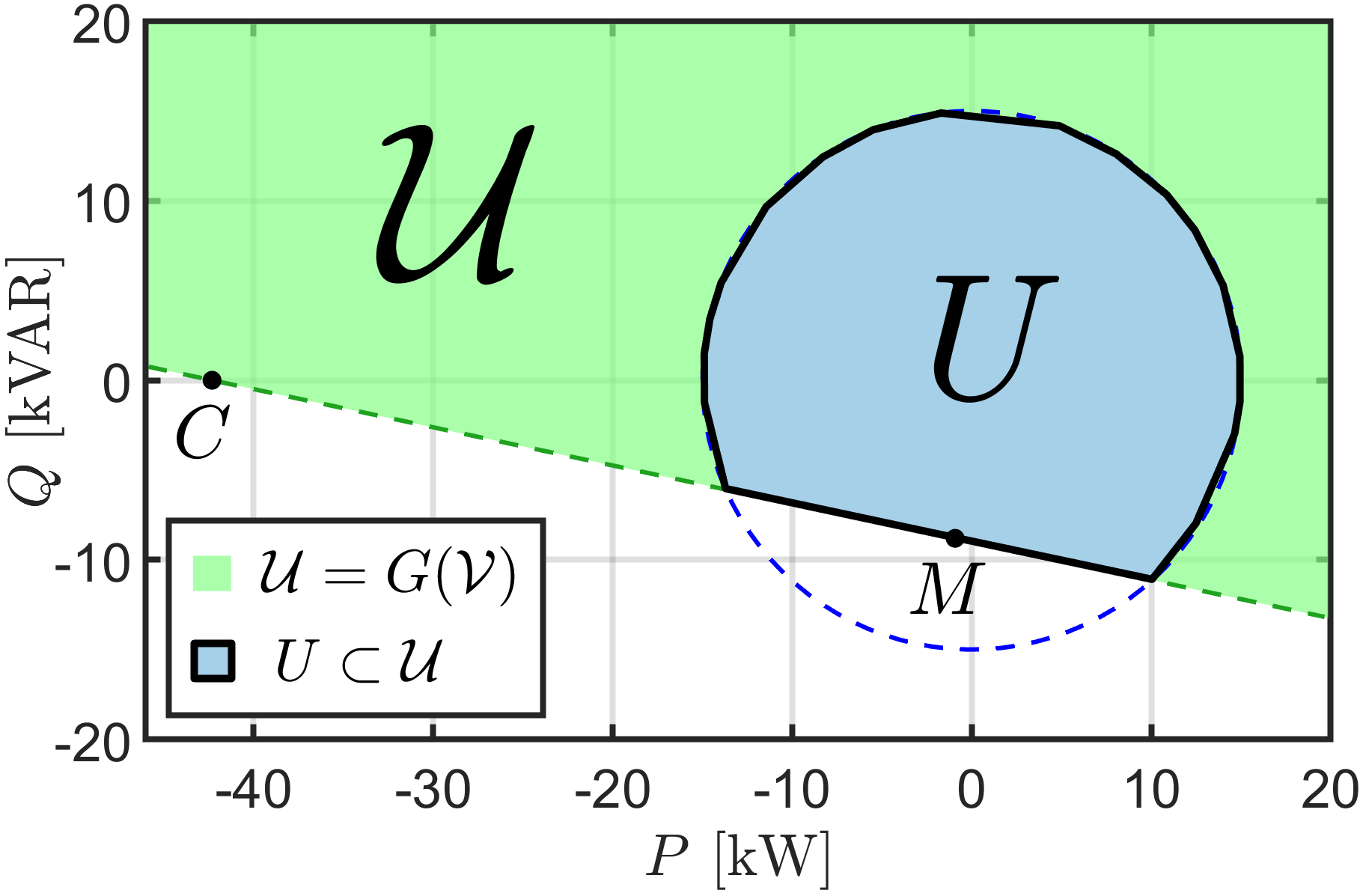}
   \caption{The set $\Uscr$ (in light green) when $\Nscr=
   G^{-1}_{\mathrm{right}}$ from \rfb{eq:N}, the set $U\subset\Uscr$
   (in light blue), and the circle $P^2+Q^2=(15\m\mathrm{kW})^2$ (in 
   dashed blue). The points $C$ and $M$ are as in \rfb{eq:CM}. (Note 
   that the lower edge of $U$ is chosen slightly above the line 
   passing through $C$ and $M$.)} \label{fig:U} \end{center}
\end{figure}

\begin{figure} \begin{center} 
   \vspace{2mm}
   \includegraphics[width=0.4\textwidth]{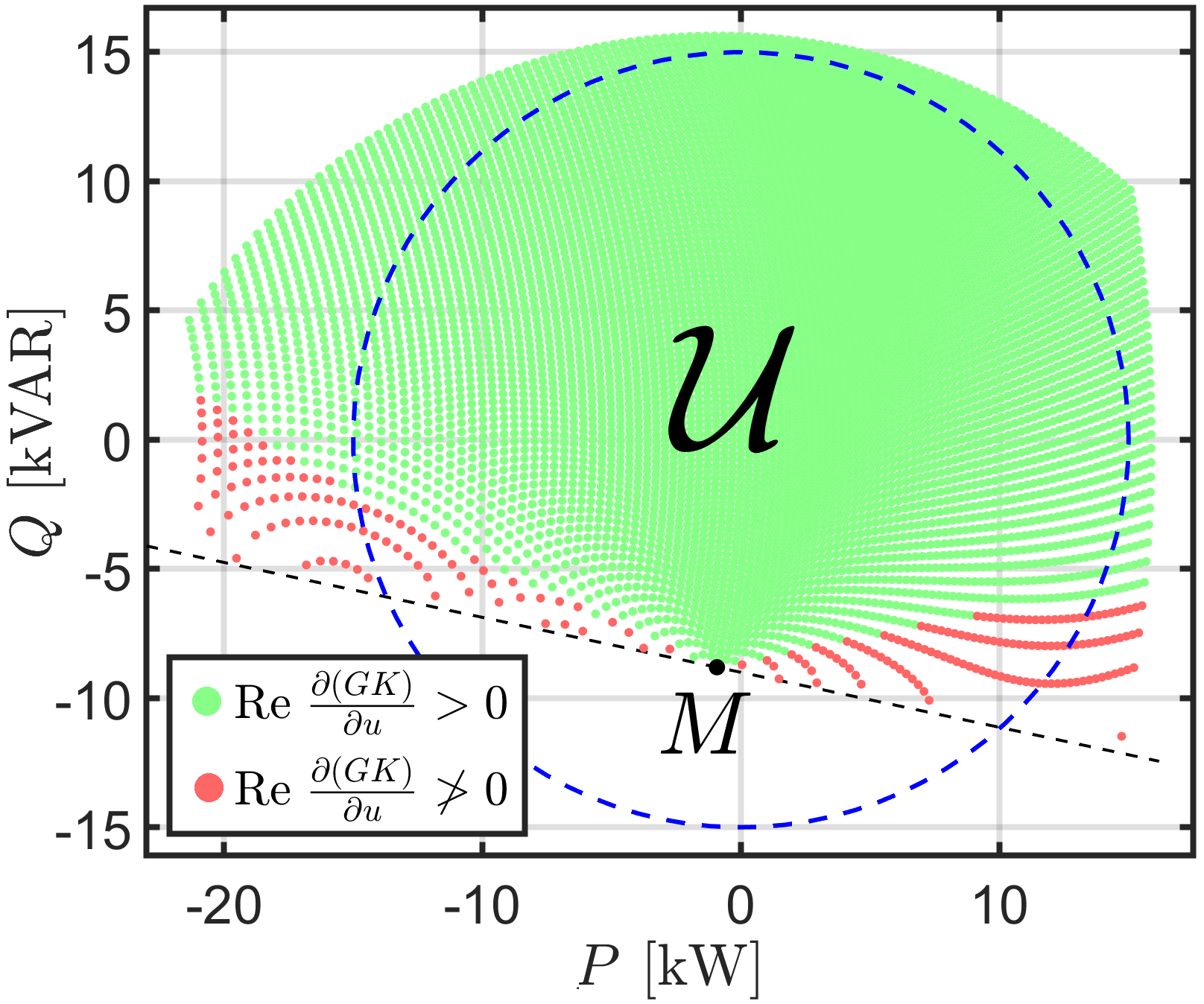}
   \caption{A glance of the set $\mathcal{U}$ (in light green) when 
   $\Nscr=K$ from \rfb{eq:K}. The circle $P^2+Q^2=(15\m\mathrm{kW})^2$ 
   is in dashed blue. Comparing the above with the set $\Uscr$ from
   Fig.~\ref{fig:U}, it is clear that a large part (depicted here in
   red) of the half-plane above the line passing through $C$ and $M$
   is missing.} \label{fig:U_K} \end{center}
\end{figure}

\begin{figure*} \centering 
   {\subfigure[The state trajectories of the saturating integrator
   $\Pi_{U}$ (in blue) and of the classical integrator (in green). We
   indicate (in red) the values of the reference $r$, to be tracked. 
   As expected, the state of the saturating integrator is never
   leaving the set $U$ from Fig.~\ref{fig:U} (shown here in light
   blue).]
   {\includegraphics[height=7.2cm]{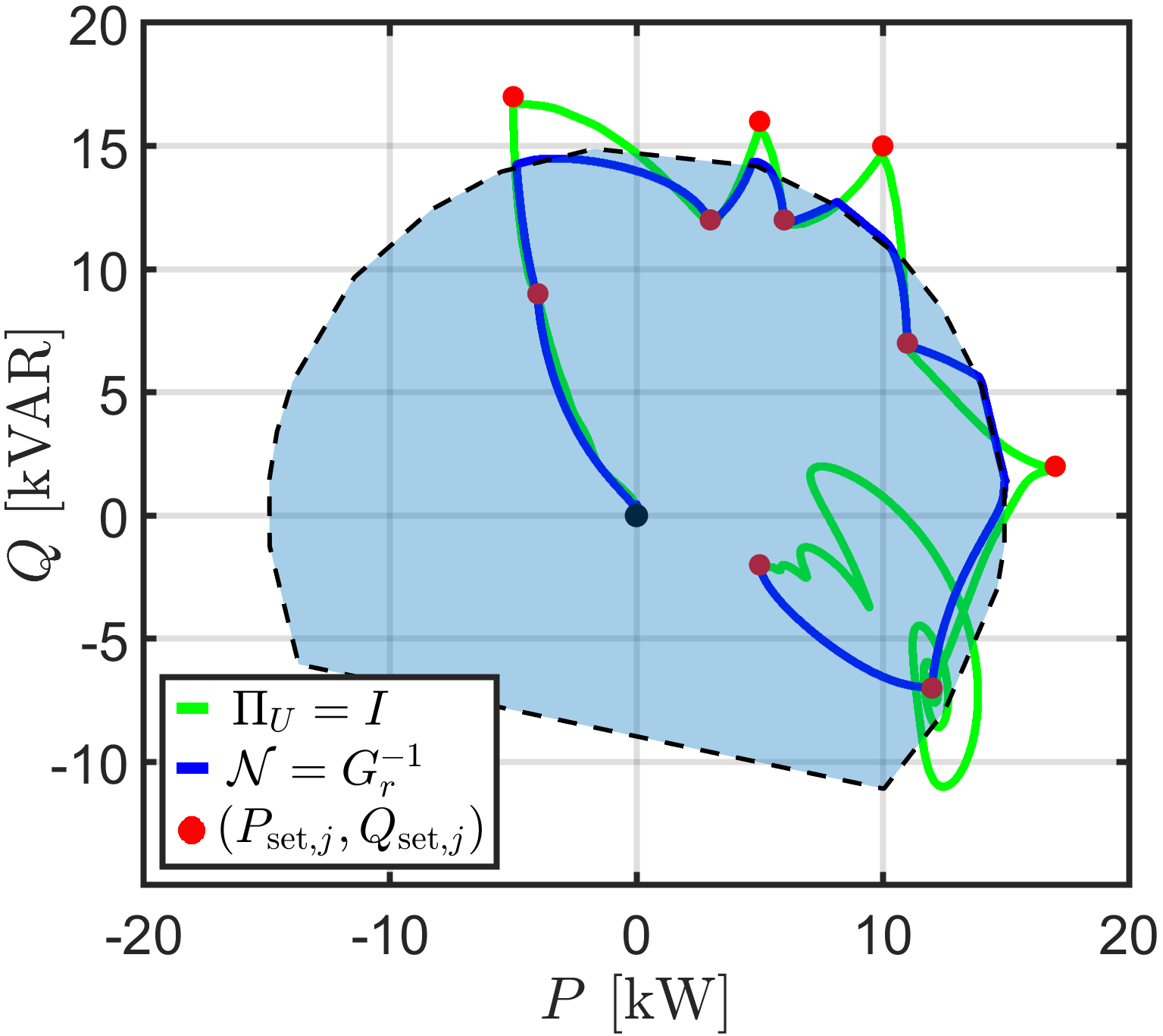}
   \label{fig:P_Q_sim}} \hfil 
   \subfigure[The values of the signal $v$ in the two closed-loop 
   systems described in the main caption. We indicate (in red) the values of
   $G^{-1}_{\mathrm{right}}(r)$, to be tracked by $v$. As expected,
   the signal $v$, in the presence of a saturating integrator, is 
   never leaving the set $V=\Nscr(U)$ (shown here in light blue).]
   {\includegraphics[height=7.2cm]{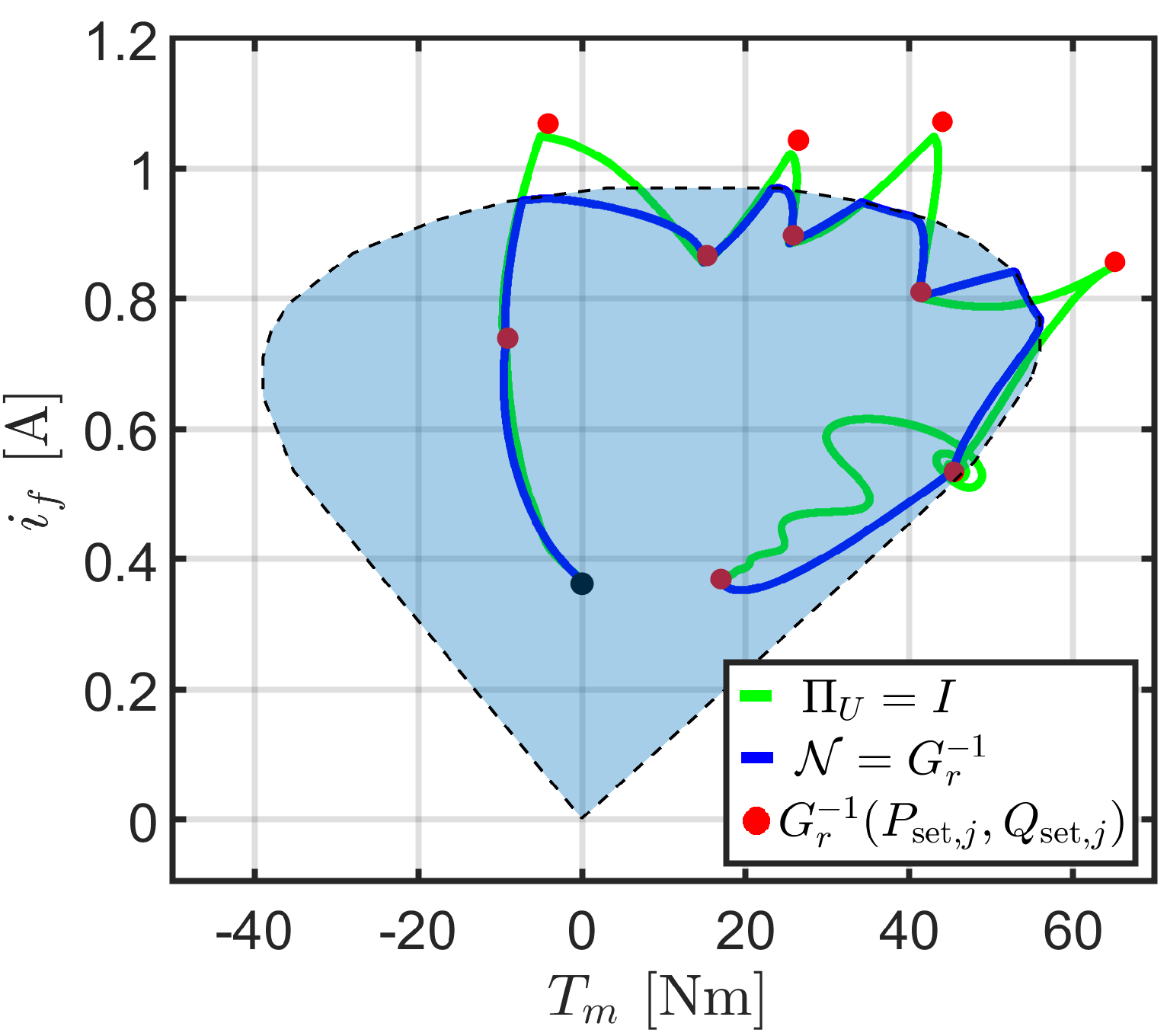}
   \label{fig:Tm_if_sim}}} 
   \caption{The comparison between the signals $u_I$ (Subfig.~a) and
   the signals $v$ (Subfig.~b) for two different closed-loop
   systems (as in Fig.~\ref{fig:clos_MIMO}): one formed by $\mathbf
   {P_0}$ in feedback with the saturating integrator $\int\Pi_U$ ($U$
   from Fig.~\ref{fig:U}), $\Nscr=G^{-1}_{\mathrm{right}}$ from
   \rfb{eq:N}, and $k=2$, whose signals are indicated in blue, and the
   other one formed by $\mathbf{P_0}$ in feedback with a classical
   integrator ($\Pi=I$), $\Nscr=K$ from \rfb{eq:K}, and $k=1$, whose
   signals are indicated in green ($\tau_p=0$ in both closed-loops).}
   \label{fig:sim}
\end{figure*}

\vspace{-2mm}
\textit{Verification of Assumption 2.} As suggested after Remark
\ref{rmk:G_inc}, we choose $\Nscr=G^{-1}_{\mathrm{right}}\in C^\infty
(\Uscr,\Vscr)$, given by \vspace{-1mm}
\begin{equation} \label{eq:N}
   \Nscr(u) \m=\m \bbm{ \ \frac{4R^2\norm{u-C}^2-V^4}{4V^2\o_g R} 
   \vspace{1mm} \ \\ \ \frac{\norm{u-M}\norm{Z}}{V\o_g m} \ }, \vspace{-3mm}
\end{equation}
where \vspace{-3mm}
\begin{equation} \label{eq:CM}
   C = \sbm{-\frac{V^2}{2R} \\ 0}, \ \ Z=\sbm{R \\ \o_g L}, \ \
   M=-\frac{V^2}{\|Z\|^2}Z,
\end{equation}
so that Assumption \ref{assumption_2} is satisfied with
$\Uscr=G(\Vscr)$. For more details on $G^{-1}_{\mathrm{right}}$, see
\cite[Theorem~3.6, Remark~3.7]{Lor2021}. The relevant portion of the set $\Uscr$ is shown 
in Fig.~\ref{fig:U}. (Outside the rectangular boundaries of Fig.~5, the powers are too large to have practical
significance for the synchronverter considered here.)
For more details on $\mathcal{U}$ see \cite[Subsect.~VI-A]{Lor2021}.

\begin{remark} \label{rmk:const_K}
We mention that an alternative choice for 
$\Nscr$ could be, e.g., $\Nscr=K\in
\rline^{2\times2}$ given by \vspace{-3mm}
\begin{equation} \label{eq:K}
   K \m=\m \bbm{\ \frac{1}{50} & 0 \ \\ \ 0 & \frac{1}{5000} \ }. \vspace{-2mm}
\end{equation}
However, this choice leads to a smaller set \m$\Uscr$, see
Fig.~\ref{fig:U_K}, where the points satisfying Assumption
\ref{assumption_2} (with $\Nscr=K$) are depicted in green. Thus,
the advantage of using $\Nscr=G^{-1}_{\mathrm{right}}$ is twofold: 
the resulting set $\Uscr$ is larger and there is no need to 
search (numerically) for the set $\Uscr$, by computing the region in 
which $\frac{\partial (G\circ K)} {\partial u}>0$, 
since $\Uscr=G(\Vscr)$.
\end{remark}

\vspace{-1mm}
\textit{The set $U$.} Due to current limitations, the safe
synchronverter operating region in the $(P,Q)$ plane is described by a
disk of radius 15\m kW. Thus, we choose $U\subset\mathcal{U}$ closed
and convex such that $P^2+Q^2\leq(15\m \mathrm{kW})^2$ in $U$. The 
set $U$ is shown (in light blue) in Fig.~\ref{fig:U}. (For 
convenience, we chose the set $U$ to be a convex polyhedron.)

\begin{table} \centering
   \renewcommand{\arraystretch}{1.5}
   \caption{Values (in chronological order) taken by the reference
   $r$.}\label{tab:r} \resizebox{\columnwidth}{!}{%
   \begin{tabular}{ccccccccccc} \hline $P_{{\rm set},\m j}$ [kW] & -4 
   & -5 & 3 & 5  & 6 & 10 & 11 & 17 & 12 & 5 \\ \hline
   $Q_{{\rm set},\m j}$ [kVAR] & 9 & 17 & 12 & 16 & 12 & 15 & 7 & 2 &
   -7 & -2 \\ \hline \end{tabular}} \vspace{1mm}
\end{table}

\subsection{Simulation results} \vspace{-2mm}

\begin{figure} \centering 
   {\subfigure[The outputs $P$ from Subfig.~\ref{fig:P_Q_sim} in time,
   	and the reference $P_{\mathrm{set},\m j}$ from Table~\ref{tab:r} (in red).]
   {\includegraphics[width=0.4\textwidth]{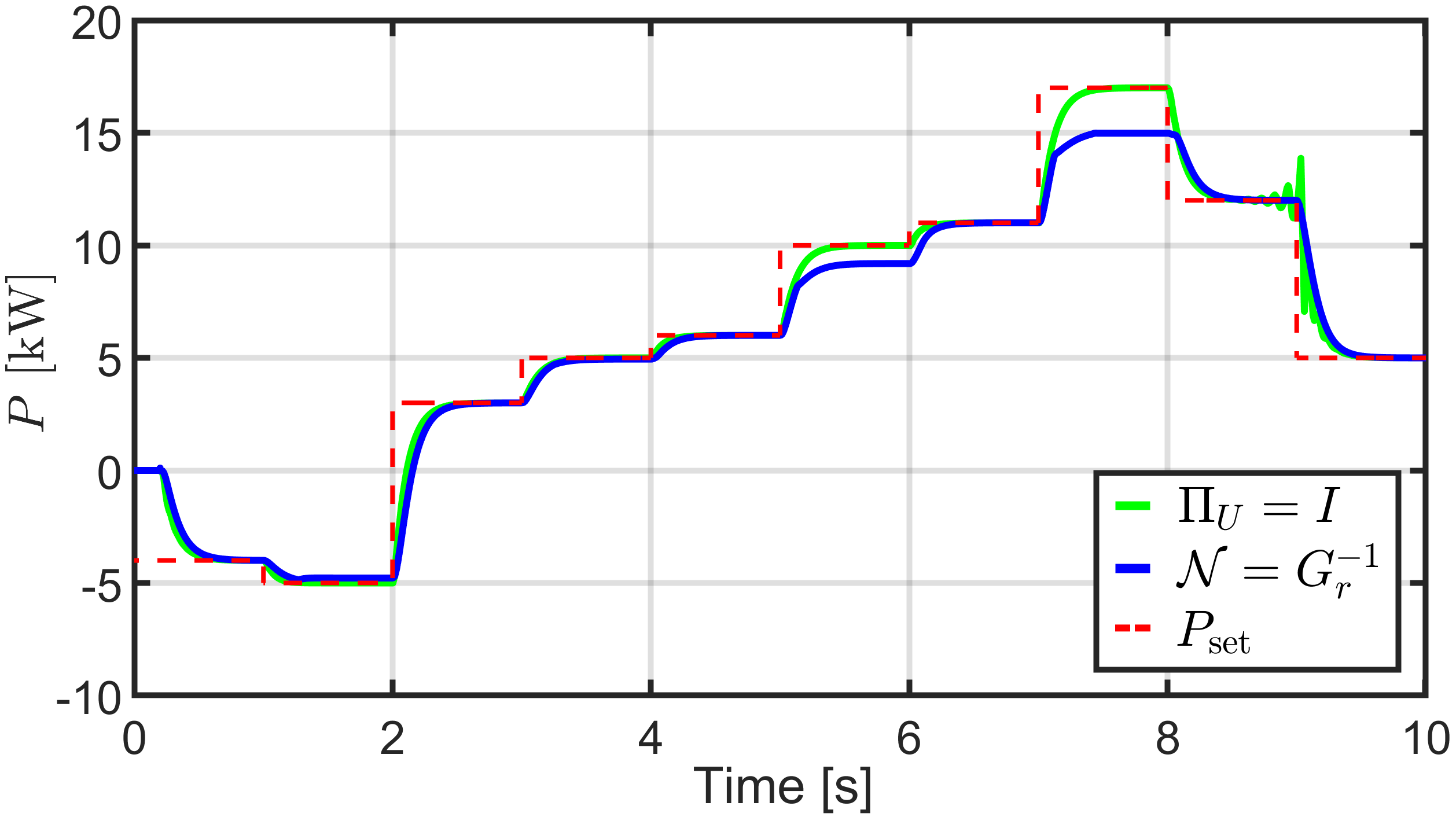}
   \label{fig:P_time}} \\
   \subfigure[The outputs $Q$ from Subfig. \ref{fig:P_Q_sim} in time,
   and the reference $Q_{\mathrm{set},\m j}$ from Table \ref{tab:r} (in red).]
   {\includegraphics[width=0.4\textwidth]{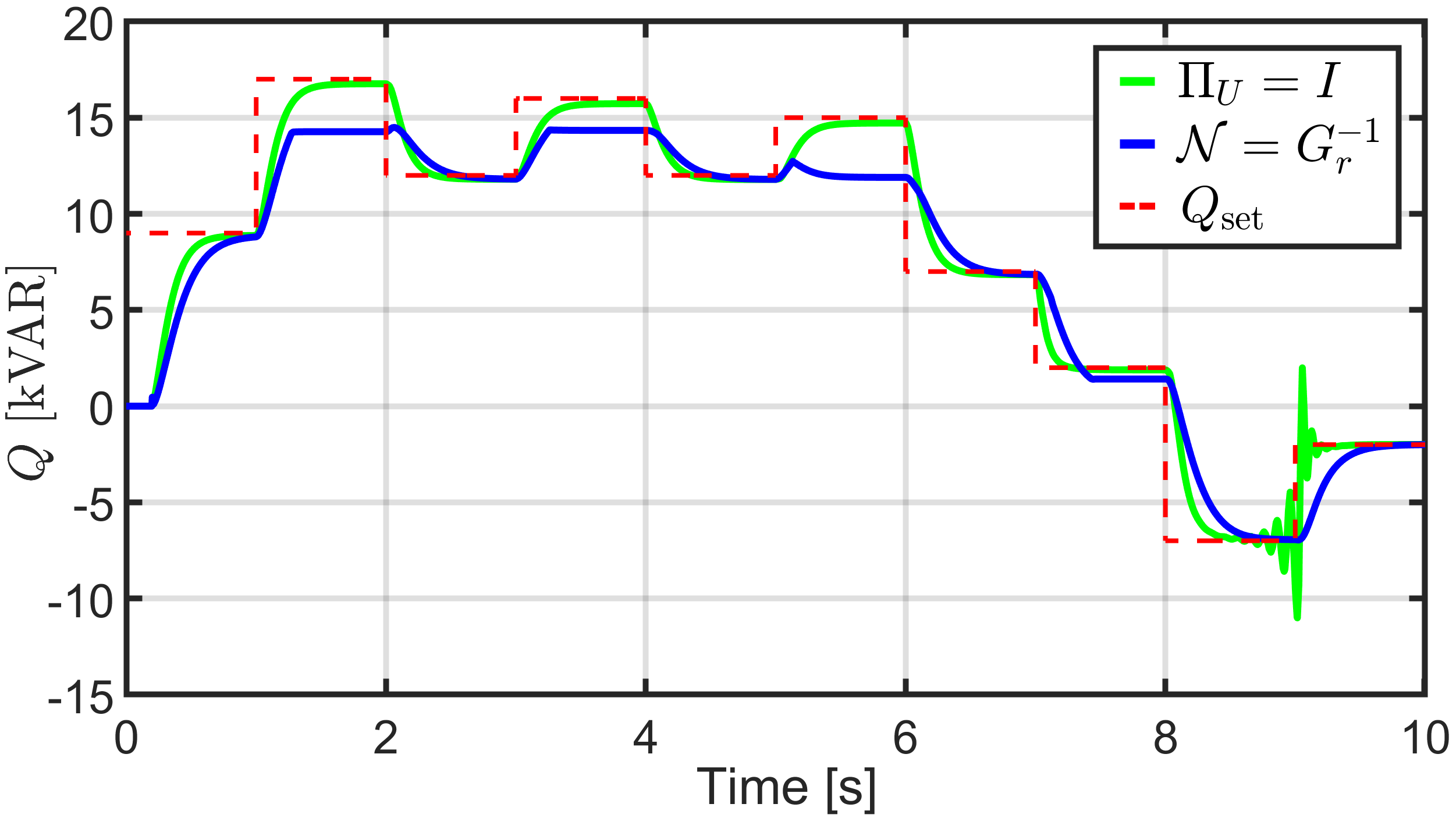}
   \label{fig:Q_time}}} \vspace{-1mm}
   \caption{The time evolution of the signals from Subfig.
   \ref{fig:P_Q_sim}.} \label{fig:time}
\end{figure}

We choose as reference signal $r$ a sequence of ten different values
for $(P_{\rm set},Q_{\rm set})$ (shown in Table~\ref{tab:r}), which
we assume to be generated by an external control loop (not modelled
here), each kept constant for 10 seconds. In Fig.~\ref{fig:P_Q_sim} we
show the comparison (in the $(P,Q)$ plane) between the state
trajectory of the saturating integrator from \rfb{eq:cl_MIMO_eq}, in
blue, with $\Nscr=G^{-1}_{\mathrm{right}}$ (given in \rfb{eq:N}) and
$k=2$, and the state trajectory of a classical integrator ($\Pi_{U}=I$ 
in \rfb{eq:cl_MIMO_eq}), in green, with $\Nscr=K$ (given in 
\rfb{eq:K}) and $k=1$ (in both cases $\tau_p=0$). It is interesting to
note that the reference point $(P_{{\rm set},\m 9},Q_{{\rm set},\m
9})=(12\m\mathrm{kW},-7\m\mathrm {kVAR})$, which generates an unstable
equilibrium point for the closed-loop system with a classical
integrator, and a stable equilibrium point for the closed-loop system
with the saturating integrator, is outside the set $\mathcal{U}$ from
Fig.~\ref{fig:U_K}, corresponding to $\Nscr=K$, but inside the set
$\Uscr$ from Fig. \ref{fig:U}, corresponding to $\Nscr=
G^{-1}_{\mathrm{right}}$. In Fig. \ref{fig:Tm_if_sim} the same
comparison is shown for the signal $v$ in the $(T_m,i_f)$ plane.
Finally, we show in Fig.~\ref{fig:P_time}, \ref{fig:Q_time}
the output active power $P$ and the output reactive power $Q$ values
(in time), for both scenarios.

\vspace{-1mm}
\begin{remark}
The step reference $r$ described above is clearly not
constant. However, it can be proved (see \cite[Prop.~4.5]{Lor2020} for the SISO case)
that the result from Theorem \ref{thm:cl_stab} can be extended for
step references (with values in $Y$) whose discontinuity
points are ``sufficiently far'' from each other.
\end{remark}

\vspace{-1mm}
\section{Conclusions} \label{sec7} 

A novel MIMO PI anti-windup controller for a stable nonlinear plant
has been proposed, based on PDS theory,
which extends our previous work \cite{Lor2020}. Under standard
assumptions, we have used SP tools to derive a
sufficient condition on the controller gain ensuring (local)
closed-loop stability and constant reference tracking. We propose
to embed the right inverse of the plant steady-state input-output map
in the controller, and we have shown the advantages of this choice
through a numerical example, namely, the output power regulation for
a grid-connected synchronverter.

\appendix
\section{Proof of Theorem \ref{thm:cl_stab}} \label{app:A}

We rewrite the closed-loop system \rfb{eq:cl_MIMO_eq} 
as a standard SP model, as in 
\cite[Sect.~III]{Lor2020}.

We introduce the variables \vspace{-1mm}
\begin{equation} \label{eq:x_tilde_u_tilde}
   \tilde{x}\vcentcolon \m=\m x-x_r, \qquad \tilde{u}_I\vcentcolon
   \m=\m u_I-u_r, \vspace{-1mm}
\end{equation}
the functions (recall $f$ from \rfb{eq:f}) \vspace{-1mm}
\begin{equation*} \label{eq:f_tilde_h_tilde}
   \begin{gathered}
   \tilde{g}(\tilde{x})\vcentcolon=g(\tilde{x}+x_r), \quad \tilde{\Pi}_{U}(\tilde{u}_I,\cdot)\vcentcolon=\m \Pi_U(\tilde{u}_I+u_r,\cdot),
   \\ \tilde{h}(\tilde{u}_I,\tilde{x})\vcentcolon=\m \tilde{\Pi}_U
   (\tilde{u}_I,r-\tilde{g}(\tilde{x})), 
   \\ \tilde{f}(\tilde{u}_I,\tilde{x})\vcentcolon=\m f(\tilde{x}+x_r,
   \tilde{u}_I+u_r), 
   \\ \tilde{\beta}(\tilde{u}_I,
   \tilde{x},k)\vcentcolon=\tilde{f}(\tilde{u}_I+\tau_pk(r-\tilde{g}
   (\tilde{x})),\tilde{x})-\tilde{f}(\tilde{u}_I,\tilde{x}),
   \end{gathered}
\end{equation*}
and we change the {\em time-scale} of \rfb{eq:cl_MIMO_eq} introducing
$s\vcentcolon= k\cdot t$. Thus, using \rfb{eq:Pi_omog}, we can rewrite
\rfb{eq:cl_MIMO_eq} as \vspace{-2mm}
\begin{equation} \label{eq:f_h_tilde_sys}
   \frac{\dd\tilde{u}_I}{\dd s} \m=\m \tilde{h}(\tilde{u}_I,
   \tilde{x}), \quad k\frac{\dd\tilde{x}}{\dd s} \m=\m \tilde{f}
   (\tilde{u}_I,\tilde{x})+\tilde{\beta}(\tilde{u}_I,\tilde{x},k).
 \vspace{-2mm}
\end{equation}
For small $k>0$, this is a standard singular perturbation model
according to \cite[Sect.~11.5]{Khalil2002}. We point out that in
\cite[Sect.~11.5]{Khalil2002} the functions describing the singularly
perturbed closed-loop systems are required to be locally Lipschitz,
which is not the case here (because of $\Pi_U$). However, our system
\rfb{eq:f_h_tilde_sys} fits the framework of
\cite[Ch.~7]{Kokotovic1999}, where it is only required that a unique
(local) closed-loop solution exists, which we have proved in
Prop.~\ref{prop:uniqueness}. \vspace{-2mm}

Following \cite[Sect.~11.5]{Khalil2002}, let \vspace{-2mm}
\begin{equation*} \label{eq:Xi_tilde}
   \tilde{\Xi}(\tilde{u}_I) \vcentcolon=\m \Xi(\Nscr
   (\tilde{u}_I+u_r))-x_r, \vspace{-2mm}
\end{equation*}
and define the fast variable \vspace{-2mm}
\begin{equation*} \label{eq:variable_change}
   \tilde{x}_f \vcentcolon=\m \tilde{x}-\tilde{\Xi}(\tilde{u}_I). 
   \vspace{-2mm}
\end{equation*}
Using the notation introduced above, we reformulate our
\rfb{eq:f_h_tilde_sys} like \cite[eqs.~(11.35),\m(11.36)]{Khalil2002},
i.e., \vspace{-3mm}
\begin{equation} \label{eq:cl_loop_u_tilde}
   \frac{\dd\tilde{u}_I}{\dd s} \m=\m \tilde{h}(\tilde{u}_I,
   \tilde{x}_f + \tilde{\Xi}(\tilde{u}_I)),
\end{equation}
\vspace{-10mm}
\begin{multline} \label{eq:cl_loop_z} 
   k\frac{\dd \tilde{x}_f}{\dd s} \m=\m \tilde{f}(\tilde{u}_I,
   \tilde{x}_f+\tilde{\Xi}(\tilde{u}_I)) + \tilde{\beta}(\tilde{u}_I,
   \tilde{x}_f+\tilde{\Xi}(\tilde{u}_I),k) \\ - k \frac{\dd
   \tilde{\Xi}}{\dd\tilde{u}_I}\tilde{h}(\tilde{u}_I,\tilde{x}_f+
   \tilde{\Xi}(\tilde{u}_I)), \vspace{-5mm}
\end{multline}
which has an equilibrium point at $(\tilde{u}_I,\tilde{x}_f)=(0,0)$.
In accordance with the change of variables \rfb{eq:x_tilde_u_tilde},
we define \vspace{-2mm}
$$ \tilde{\mathcal{U}}\vcentcolon=\m \mathcal{U}-u_r\subset\rline^p
\quad \text{and} \quad \tilde{U}\vcentcolon=\m 
U-u_r\subset\tilde{\mathcal{U}},\vspace{-2mm}$$
which contain the origin. Thus, the state space of the closed-loop
system \rfb{eq:cl_loop_u_tilde},\rfb{eq:cl_loop_z} is 
$\tilde{\mathcal{X}}\vcentcolon=\tilde{\mathcal{U}}\times\rline^n$.

\vspace{-1mm}
Using standard arguments, see \cite[Ch.~11]{Khalil2002},
\cite[Ch.~7]{Kokotovic1999} or \cite[Sect.~III]{Lor2020}, we
identify the reduced model and the boundary-layer system associated to
\rfb{eq:cl_loop_u_tilde}-\rfb{eq:cl_loop_z}.  Recall $G$ from
Assumption~\ref{assumption_2}. Define the function \vspace{-2mm}
\begin{equation*} \label{eq:S_G_tilde}
   \tilde{G}(\tilde{u}_I)\vcentcolon = \tilde{g}(\tilde{\Xi}
   (\tilde{u}_I)) = G(\Nscr(\tilde{u}_I+u_r)). \vspace{-2mm}
\end{equation*}
The \textit{reduced (slow) model} associated to 
\rfb{eq:cl_loop_u_tilde}-\rfb{eq:cl_loop_z} is
obtained by taking $\tilde{x}_f=0$ in 
\rfb{eq:cl_loop_u_tilde}, which leads to \vspace{-2mm}
\begin{equation} \label{eq:reduced_model}
   \frac{\dd\tilde{u}_I}{\dd s} \m=\m \tilde{\Pi}_{U}(\tilde{u}_I,r-
\tilde{G}(\tilde{u}_I)). \vspace{-5mm}
\end{equation}

The \textit{boundary-layer (fast) system} associated to 
\rfb{eq:cl_loop_u_tilde}-\rfb{eq:cl_loop_z} is obtained by
rewriting \rfb{eq:cl_loop_z} in the original fast time scale $t$
and then taking $k=0$, which yields  \vspace{-1mm}
\begin{equation} \label{eq:fast_system}
   \dot{\tilde{x}}_f\m=\m \tilde{f}(\tilde{u}_I,\tilde{x}_f +
   \tilde{\Xi}(\tilde{u}_I)), \vspace{-1mm}
\end{equation} 
where $\tilde{u}_I\in\tilde{\Uscr}$ is treated as a fixed parameter.


\vspace{-1mm}
We are now ready to prove the stability of the equilibrium point
$(x_r,u_r)$ of the closed-loop system \rfb{eq:cl_MIMO_eq}, using 
SP theory. We follow the
arguments in \cite[Sect.~IV]{Lor2020}, which are based on the
guidelines of \cite[Sect.~11.5]{Khalil2002}. (Note that the fast
variable $z$ in \cite{Lor2020} is denoted here by $\tilde{x}_f$.)

Define the set $\tilde{U}_\delta\vcentcolon=\tilde{U}+B_\delta$, where
$B_\delta$ denotes the closed ball of radius $\delta>0$ in $\rline^p$.
We choose $\delta$ such that $\tilde{U}_\delta\subset\tilde{\Uscr}$.
We will use \cite[Th.~4.2]{Lor2020}, but with $\tilde{U}_\delta\subset
\rline^p$ (instead of $\tilde{U}_\delta\subset\rline$). To check this
extension of \cite[Th.~4.2]{Lor2020}, it is enough to replace the 
Lipschitz property of the saturating integrator $\mathscr{S}$ in 
\cite{Lor2020}, with the contraction property \rfb{eq:contraction} of
the operator $\Pi_U$ in the proof of \cite[Th.~4.2]{Lor2020}.

{\em Step 1: Stability of the reduced model \rfb{eq:reduced_model}.}
Let $r\in Y$. The aforementioned extension of \cite[Theorem~4.2]
{Lor2020} demands the existence of a Lyapunov function $V$ for
\rfb{eq:reduced_model} (defined on $\tilde{U}_\delta$) such that
\vspace{-2mm}
\begin{equation} \label{eq:V_condit}
   \begin{gathered} c_1\norm{\tilde{u}_I}^2\leq V(\tilde{u}_I)
   \leq c_2 \norm{\tilde{u}_I}^2, \\ \frac{\dd V}{\dd \tilde{u}_I}
   \tilde{h}(\tilde{u}_I,\tilde{\Xi}(\tilde{u}_I))\leq-c_3\norm{
   \tilde{u}_I}^2,\ \ \ \bigg\|\frac{\dd V}{\dd \tilde{u}_I}\bigg\|
   \leq c_4 \norm{\tilde{u}_I},\end{gathered}
\end{equation}
for all $\tilde{u}_I\in\tilde{U}_\delta$, where $c_1,\dots c_4$ are
positive constants.

\vspace{-2mm}
As in \cite[Subsec.~IV-A]{Lor2020}, we consider the candidate Lyapunov
function \vspace{-2mm}
\begin{equation*} \label{eq:Lyap_red_mod}
   V(\tilde{u}_I)=\frac{1}{2}\norm{\tilde{u}_I}^2 \FORALL \tilde{u}_I
   \in\tilde{U}_\delta. \vspace{-2mm}
\end{equation*}
Its derivative along the trajectories of \rfb{eq:reduced_model} is 
\vspace{-2mm}
\begin{equation*}
   \frac{\dd V}{\dd s}=\langle\m\tilde{\Pi}_{U}(\tilde{u}_I,\tilde{G}
   (0)-\tilde{G}(\tilde{u}_I)),\tilde{u}_I\m\rangle. \vspace{-2mm}
\end{equation*}
The (unique) equilibrium point of \rfb{eq:reduced_model} is $0\in {\rm
int}\m\tilde{U}_\delta$ and $\tilde{G}$ is strictly monotone from
Assumption \ref{assumption_2}. Therefore, the operator $\tilde{\Pi}_U$
behaves like the identity and the block $\int\tilde{\Pi}_U$ reduces to
a classical integrator. Thus \vspace{-2mm}
$$ \frac{\dd V}{\dd s} \m=\m \langle\m \tilde{G}(0) - \tilde{G}
   (\tilde{u}_I),\tilde{u}_I\m\rangle\leq-\mu\norm{\tilde{u}_I}^2,
   \vspace{-2mm}$$
and the conditions \rfb{eq:V_condit} are easily seen to hold.

\vspace{-1mm}
\textit{Step 2: Stability of the boundary-layer system
\rfb{eq:fast_system}.} \cite[Theorem~4.2]{Lor2020} requires the
existence of a Lyapunov function $W$ for \rfb{eq:fast_system} (defined
on $\tilde{U}_\delta\times B_{\e_0}$) such that \vspace{-2mm}
\begin{equation} \label{eq:W_lemma_9_8}
   \begin{gathered} b_1\norm{\tilde{x}_f}^2\leq W(\tilde{u}_I,
   \tilde{x}_f)\leq b_2 \norm{\tilde{x}_f}^2, \\ \frac{\partial W}
   {\partial \tilde{x}_f} \tilde{f}(\tilde{u}_I,\tilde{x}_f+
   \tilde{\Xi}(\tilde{u}_I)) \leq -b_3\norm{\tilde{x}_f}^2, \\
   \bigg\|\frac{\partial W}{\partial \tilde{x}_f}\bigg\| \leq b_4
   \norm{\tilde{x}_f},\quad \bigg\|\frac{\partial W}{\partial 
   \tilde{u}_I}\bigg\| \leq b_5 \norm{\tilde{x}_f}^2, \vspace{-1mm}
   \end{gathered}
\end{equation}
for all $(\tilde{u}_I,\tilde{x}_f)\in\tilde{U}_\delta\times B_{\e_0}$
(recall $\e_0$ from Assumption \ref{assumption_1}), where $b_1,\dots
b_5$ are positive constants. As in \cite{Lor2020}, we want to use
\cite[Lemma~9.8]{Khalil2002} to guarantee the existence of a function
$W$ such that \rfb{eq:W_lemma_9_8} holds. To check its assumptions, we
use the arguments of \cite[Subsec.~IV-B]{Lor2020}, with the difference
that here $\tilde{U}_\delta\subset\rline^p$ (instead of
$\tilde{U}_\delta\subset\rline$). Thus, we can simply replace
$F_j(z,\tilde{u}_I)$ there with $F_{lj}(\tilde{x}_f,\tilde{u}_I)
\vcentcolon=\frac{\partial p_j}{\partial \tilde{u}_{I_l}}$ here, for
all $l\in\{1,2,\dots \m p\}$ and for all $j\in\{1,2,\dots \m n\}$, to
guarantee that the assumptions of \cite[Lemma~9.8]{Khalil2002} are
met and, thus, that a function $W$ satisfying
\rfb{eq:W_lemma_9_8} exists.

\vspace{-1mm}
\textit{Step 3: Stability of the closed-loop system
\rfb{eq:cl_MIMO_eq}.} We complete the proof of Theorem
\ref{thm:cl_stab} by following step-by-step that of
\cite[Th.~4.3]{Lor2020}, and using the extension of \cite[Th.~4.2]{Lor2020} discussed before
Step 1. \hfill $\blacksquare$


\small{ {\bf Pietro Lorenzetti} received the MEng degree in Mechatronic
Eng. from Politecnico di Torino, and in Automation and Control Eng.
from Politecnico di Milano in 2017,
with honours, thanks to the double-degree program ``Alta Scuola Politecnica''.
He is currently an Early Stage Researcher within the Marie Curie ITN project ``ConFlex''
in Tel Aviv University,
under the supervision of G.~Weiss. He is the recipient of the 
IFAC Young Author Award of the MICNON2021 conference.
His research interests include
nonlinear systems, nonlinear control, and power systems stability.}

\small{ {\bf George Weiss} received the MEng degree in control
engineering from the Polytechnic Institute of Bucharest, Romania, in
1981, and the Ph.D. degree in applied mathematics from the Weizmann
Institute, Rehovot, Israel, in 1989. He was with Brown University,
Providence, RI, Virginia Tech, Blacksburg, VA, Ben-Gurion University,
Beer Sheva, Israel, the University of Exeter, U.K., and Imperial
College London, U.K. His current research interests include
distributed parameter systems, operator semigroups, passive and
conservative systems (linear and nonlinear), power electronics,
microgrids, repetitive control, sampled data systems, and wind-driven
power generators. He is leading research projects for the European
Commission and for the Israeli Ministry of Infrastructure, Energy and
Water.}


\begin{thebibliography}{xx} \vspace{-2mm}

\bibitem[{\AA}str{\"o}m and Rundqwist(1989)]{Astrom1989} 
 K.~J.~{\AA}str{\"o}m, and L.~Rundqwist.
\newblock ``Integrator windup and how to avoid it,''
\newblock \emph{1989 ACC. IEEE}, pp.~1693-1698, 1989.


\bibitem[Davison(1976)]{Davison1976} E.~J.~Davison.
\newblock ``Multivariable tuning regulators: The feedforward and
 robust control of a general servomechanism problem,''
\newblock \emph{IEEE TAC}, vol.~21, no.~1,
 pp.~35-47, 1976.

\bibitem[Desoer and Lin(1985)]{Desoer1985} C.~Desoer and C.A.~Lin.
\newblock ``Tracking and disturbance rejection of MIMO nonlinear 
 systems with PI controller'',
\newblock \emph{IEEE TAC}, vol.~30, no.~9, 
 pp.~861--867, 1985.


\bibitem[Edwards and Postlethwaite(1998)]{Edwards1998} C.~Edwards, 
 and I.~Postlethwaite.
\newblock ``Anti-windup and bumpless-transfer schemes,''
\newblock \emph{Automatica}, vol.~34, pp.~199-210, 1998.

\bibitem[Francis(1975)]{Francis1975} B.~A.~Francis and W.~M.~Wonham.
\newblock ``The internal model principle for linear multivariable
 regulators,'' 
\newblock {\em Appl. Math. Optim.}, vol.~2, pp.~170-194, 1975.

\bibitem[Guiver et~al.(2017) Guiver, Logemann, and 
 Townley]{Guiver2017} C.~Guiver, H.~Logemann, and S.~Townley.
\newblock ``Low-gain integral control for multi-input multi-output
 linear systems with input nonlinearities'',
\newblock \emph{IEEE TAC}, vol.~62, 
 pp.~4776-4783, 2017.

\bibitem[Hauswirth et~al.(2021) Hauswirth, Bolognani, and 
 D\"{o}rfler]{Hauswirth2021} A.~Hauswirth, S.~Bolognani, and 
 F.~D\"{o}rfler.
\newblock ``Projected dynamical systems on irregular, non-Euclidean
 domains for nonlinear optimization,'' 
\newblock {\em SIAM Journal on Control and Optimization}, vol.~59,
pp.~635-668, 2021.



\bibitem[Hauswirth et~al.(2020) Hauswirth, D\"{o}rfler, and 
 Teel]{Hauswirth2020b} A.~Hauswirth, F.~D\"{o}rfler, and A.~Teel.
\newblock ``On the robust implementation of projected dynamical 
 systems with anti-windup controllers,''
\newblock {\em Proc. of the 2020 ACC}, pp.~1286-1291, 2020.

\bibitem[Huang et al.~(2019)]{Huang2019} X.~Huang, H.~K.~Khalil and Y.~Song, 
\newblock ``Regulation of nonminimum-phase nonlinear systems using slow integrators and high-gain feedback,''
\newblock \emph{IEEE Trans. on Automatic Control}, vol.~64, pp.~640-653, 2019.


\bibitem[Khalil(2002)]{Khalil2002} H.~K.~Khalil.
\newblock {\em Nonlinear Systems; 3rd ed.}
\newblock Prentice-Hall, Upper Saddle River, NJ, 2002.

\bibitem[Kothare et~al.(1994)]{Kothare1994}
 M.~V.~Kothare, P.J.~Campo, M.~Morari, and C.N.~Nett.
\newblock ``A unified framework for the study of anti-windup
 designs,''
\newblock \emph{Automatica}, vol.~30, pp.~1869-1883, 1994. 

\bibitem[{Kokotović et~al.(1999) 
 Kokotović, Khalil, and O'Reilly}]{Kokotovic1999} 
 P.~Kokotović, H.~K.~Khalil, and J.~O'Reilly.
\newblock {\em Singular Perturbation Methods in Control: Analysis 
 and Design}.
\newblock SIAM, 1999.

\bibitem[Konstantopoulos et~al.(2016) Konstantopoulos, Zhong, Ren, 
 Krstic]{Konstantopoulos2016} G.~C.~Konstantopoulos, Q.-C.~Zhong, 
 B.~Ren, and M.~Krstic.
\newblock ``Bounded integral control of input-to-state practically
 stable nonlinear systems to guarantee closed-loop stability'',
\newblock \emph{IEEE TAC}, vol.~61, 
 pp.~4196-4202, 2016.


\bibitem[Logemann et~al~(1999)]{Logemann1999} 
 H.~Logemann, E.P.~Ryan and S.~Townley.
\newblock ``Integral control of linear systems with actuator
 nonlinearities: lower bounds for the maximal regulating gain'',
\newblock \emph{IEEE TAC}, vol.~44, 
 pp.~1315-1319, 1999.
 
\bibitem[Lorenzetti et~al(2022) Lorenzetti, Kustanovich, Shivratri,
 and Weiss]{Lor2021} P.~Lorenzetti, Z.~Kustanovich, S.~Shivratri, and G.~Weiss.
\newblock ``The equilibrium points and stability of grid-connected
 synchronverters,''
\newblock {\em IEEE Trans. Power Systems}, vol.~37, pp.~1184-1197, 2022.

\bibitem[Lorenzetti and Weiss(2022)]{Lor2020} 
 P.~Lorenzetti and G.~Weiss.
\newblock ``Saturating PI control of stable nonlinear systems using
 singular perturbations,''
\newblock {\em IEEE TAC}, early access, 2022.

\bibitem[Lorenzetti and Weiss(2021)]{Lor2021Conf}
 P.~Lorenzetti and G.~Weiss.
\newblock ``Integral control of stable MIMO nonlinear systems with
 input constraints,''
\newblock to appear in the \textit{Proc. of the $3^{rd}$ MICNON
 Conference, Tokyo, September, 2021}.

\bibitem[Lorenzetti et~al(2020) Lorenzetti, Weiss, and
 Natarajan]{Lor2020Conf} 
 P.~Lorenzetti, G.~Weiss and V.~Natarajan.
\newblock ``Integral control of stable nonlinear systems based on
 singular perturbations'', 
\newblock \textit{IFAC-PapersOnLine}, vol.~53, pp.~6157-6164, 2020.


\bibitem[Morari(1985)]{Morari1985} M.~Morari.
\newblock ``Robust stability of systems with integral control,''
\newblock \emph{IEEE TAC}, vol.~30, 
 pp.~574-577, 1985.

\bibitem[Nagurney and Zhang(1995)]{Nagurney1995} 
 A.~Nagurney and D.~Zhang. 
\newblock \emph{Projected Dynamical Systems and Variational 
 Inequalities with Applications},
\newblock Springer Science \& Business Media, 1995.

\bibitem[Natarajan and Weiss(2018)]{Natarajan2018}
 V.~Natarajan and G.~Weiss.
\newblock ``Almost global asymptotic stability of a grid-connected
 synchronous generator,''
\newblock {\em Math. of Control, Signals and Systems}, vol.~30, 
 2018.

\bibitem[Natarajan and Weiss(2017)]{Natarajan2017} 
 V.~Natarajan, and G.~Weiss.
\newblock ``Synchronverters with better stability due to virtual 
 inductors, virtual capacitors, and anti-windup'', 
\newblock \emph{IEEE Trans. on Industrial Electronics}, vol.~64, 
 pp.~5994-6004, 2017.


\bibitem[Simpson-Porco(2020)]{Simpson-Porco2020} J.~W.~Simpson-Porco.
\newblock ``Analysis and synthesis of low-gain integral controllers
 for nonlinear systems,''
\newblock \emph{IEEE TAC}, published 
 online, 2020.

\bibitem[Simpson-Porco(2021)]{Simpson-Porco2021} J.~W.~Simpson-Porco.
\newblock ``Low-gain stability of projected integral control for 
 input-constrained discrete-time nonlinear systems,'' 
\newblock {\em IEEE Control Systems Letters}, vol.~6, pp.~788-793,
 2021.



\bibitem[Sontag and Wang(1996)]{Sontag1996}
E.~D.~Sontag, and Y.~Wang.
\newblock ``New characterizations of input to state stability'',
\newblock \emph{IEEE TAC},
vol.~41, pp.~1283-1294., 1996.

\bibitem[Tarbouriech and Turner(2009)]{Tarbouriech2009} 
 S.~Tarbouriech, and M.~Turner.
\newblock ``Anti-windup design: an overview of some recent advances
 and open problems,''
\newblock \emph{IET Control Theory \& Appl.}, vol.~3, pp.~1-19, 2009.


\bibitem[Teo and How(2011)]{Teo2011} J.~Teo and J.P.~How.
\newblock ``Region of attraction comparison for gradient projection 
 anti-windup compensated systems,''
\newblock {\em 2011 50th IEEE CDC}, pp.~5509-5515, 2011.


\bibitem[Wang et~al.(2020) Wang, Ren, Zhong, and Dai]{Wang2020}
 Y.~Wang, B.~Ren, Q.-C.~Zhong, and J.~Dai.
\newblock ``Bounded integral controller with limited control power
 for nonlinear multiple-input multiple-output systems,'' 
\newblock {\em IEEE Trans. on Cont. Systems Tech.}, early 
 access, 2020.

\bibitem[Zaccarian and Teel(2002)]{Zaccarian2002}
 L.~Zaccarian, and A.~R.~Teel.
\newblock ``A common framework for anti-windup, bumpless transfer
 and reliable designs'',
\newblock \emph{Automatica}, vol.~38, pp.~1735-1744, 2002.


\bibitem[Zhong and Weiss (2011)]{Zhong2011}
 Q.-C.~Zhong and G.~Weiss.
\newblock ``Synchronverters: Inverters that mimic synchronous
 generators,'' 
\newblock {\em IEEE Trans. Industr. Electronics}, vol.~58, 
 pp.~1259-1267, 2011.

\end{thebibliography}
\end{document}